\documentclass[a4wide, oneside]{article}  

\usepackage{myPackages}
\usepackage{myAssignments3}
\usepackage{myOperators}

\renewcommand{\S}{\Sigma}
\renewcommand{\t}{\tau}

\usepackage{tikz}
\usetikzlibrary{arrows,decorations.pathmorphing,backgrounds,positioning,fit,petri}
\usetikzlibrary{calc,arrows.meta,positioning}
\usepackage{ragged2e}

\usepackage{fancyhdr,floatpag}
\floatpagestyle{fancy}% Page style for float-page only

\usepackage{caption}
\usepackage{subcaption}

\titleformat{\section}
{\normalfont\normalsize\bfseries\center}{\thesection}{1em}{}

\titleformat{\subsection}
{\normalfont\normalsize\bfseries}{\thesection}{1em}{}

\begin{document}
\renewcommand{\abstractname}{\vspace{-\baselineskip}}

\title{Contact surgery on the Hopf link: classification of fillings}

\author{Edoardo Fossati	}

\date{}

\newcommand{\Addresses}{{% additional braces for segregating \footnotesize
  \bigskip
 % \footnotesize

  E.~Fossati, \textsc{Scuola Normale Superiore, Piazza dei Cavalieri 7, Pisa.}\par\nopagebreak
  \textit{E-mail address:} \texttt{edoardo.fossati@sns.it}

}}

\maketitle

\begin{abstract}
Let $H\subseteq S^3$ be the two-component Hopf link. After choosing a Legendrian representative of $H$ with respect to the standard tight contact structure on $S^3$ we perform contact $(-1)$-surgery on the link itself. The manifold we get is a lens space together with a tight contact structure on it, which depends on the chosen Legendrian representative. We classify its minimal symplectic fillings up to homeomorphism (and often up to diffeomoprhism), extending the results of \cite{lisca} which covers the case of universally tight structures and the article of \cite{plamenVHM} which describes the fillings of $L(p,1)$.
\end{abstract}

\section{Introduction}

The problem of classifying Stein (or more generally, symplectic) fillings of contact 3-manifolds has come to the attention of topologists since the pioneering work of Elisahberg \cite{eliashberg}. Showing that $S^3$ with its standard tight contact structure has a unique Stein filling, which is $(D^4,J_{std})$, is the starting point of this active research area. In the last years several other works appeared on this subject, including McDuff's article \cite{mcduff} that classifies the fillings of $(L(p,1),\xi_{std})$, later extended by Lisca in \cite{lisca} to any $(L(p,q),\xi_{std})$. Trying to get a complete list of the possible tight contact structures (up to isotopy) on a 3-manifold is in general very hard, but at least on lens spaces this list is available thanks to the work of Honda \cite{honda}. This gives extra motivation to study these spaces, in order to fully understand the connection to the 4-dimensional symplectic world. A key distinction within tight contact structures is the one between \emph{universally tight} and \emph{virtually overtwisted} structures, see definition \ref{utvot}. In the articles of McDuff and Lisca quoted above the \emph{standard} structure refers to the universally tight one(s), but Honda showed that there are many more structures, which are indeed virtually overtwisted.

\paragraph{Acknowledgments} The author wishes to thank Paolo Lisca for introducing him to this topic and for the support during the development of this work, which will be part of his PhD thesis. Many thanks are due to Andr\'as Stipsicz for all the useful corrections and comments on the first drafts of this article. 

\paragraph{Known results} Plamenevskaya and Van Horn Morris \cite{plamenVHM} showed that the virtually overtwisted structures on $L(p,1)$ have a unique Stein filling. Another classification result about fillings of virtually overtwisted structures on certain families of lens spaces is given by Kaloti \cite{kaloti}. The goal of this article is to prove the following:

\begin{thm}\label{theorem1} Let $L$ be the lens space resulting from Dehn surgery on the Hopf link with framing $-a_1$ and $-a_2$, with $a_1,a_2\geq 2$. Let $\xi_{vo}$ be a virtually overtwisted contact structure on $L$. 
Then $(L,\xi_{vo})$ has:
\begin{itemize}
\item a unique (up to diffeomorpism) Stein filling if $a_1\neq 4\neq a_2$,
\item two homeomorphism classes of Stein fillings, distinguished by the second Betti number $b_2$, if at least one of $a_1$ and $a_2$ is equal to 4 and the corresponding rotation number is $\pm 2$. Moreover, the diffeomorphism type of the Stein filling with bigger $b_2$ is unique. If the rotation number is not $\pm 2$, then we have again a unique filling.
\end{itemize}
\end{thm}

\begin{figure}[ht!]
\centering
\includegraphics[scale=0.4]{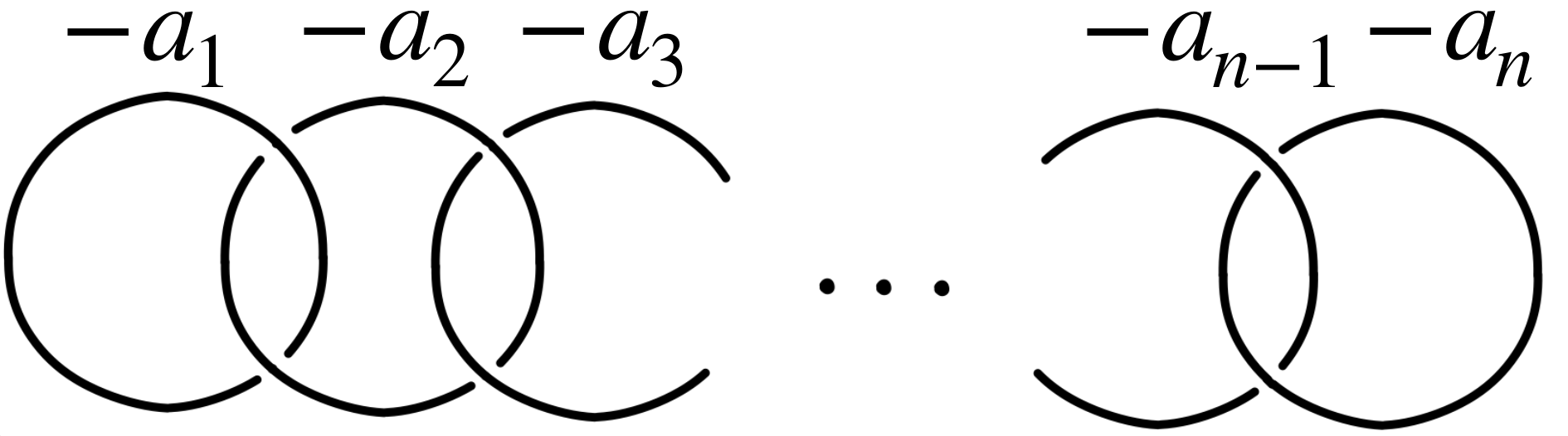}
\caption{Surgery link producing $L(p,q)$.}
\label{chain}
\end{figure}

\noindent Given a pair of coprime integers $p>q>1$, we can consider the continued fraction expansion
\[\frac{p}{q}=[a_1,a_2,\ldots a_n]=a_1-\frac{1}{a_2-\frac{1}{\ddots -\frac{1}{a_n}}},\]
with $a_i\geq 2$ for every $i$. 
As a smooth oriented 3-manifold, $L(p,q)$ is the integral surgery on a chain of unknots with framings $-a_1,-a_2,\ldots -a_n$ (figure \ref{chain}). 

\noindent Before equipping $L(p,q)$ with a tight contact structure, we recall a further dichotomy in contact topology:

\begin{defn}\label{utvot} A tight contact strcutre $\xi$ on a 3-manifold $M$ is called \emph{universally tight} if its pullback to the universal cover $\widetilde{M}$ is tight. The tight structure $\xi$ is called \emph{virtually overtwisted} if its pullback to some finite cover is overtwisted.
\end{defn}

\begin{nrem}
A consequence of the geometrization conjecture is that the fundamental group of any 3-manifold is \emph{residually finite} (i.e. any non trivial element is in the complement of a normal subgroup of finite index), and this implies that any tight contact structure is either universally tight or virtually overtwisted (see \cite{honda}).
\end{nrem}

\begin{figure}[ht!]
\centering
\includegraphics[scale=0.4]{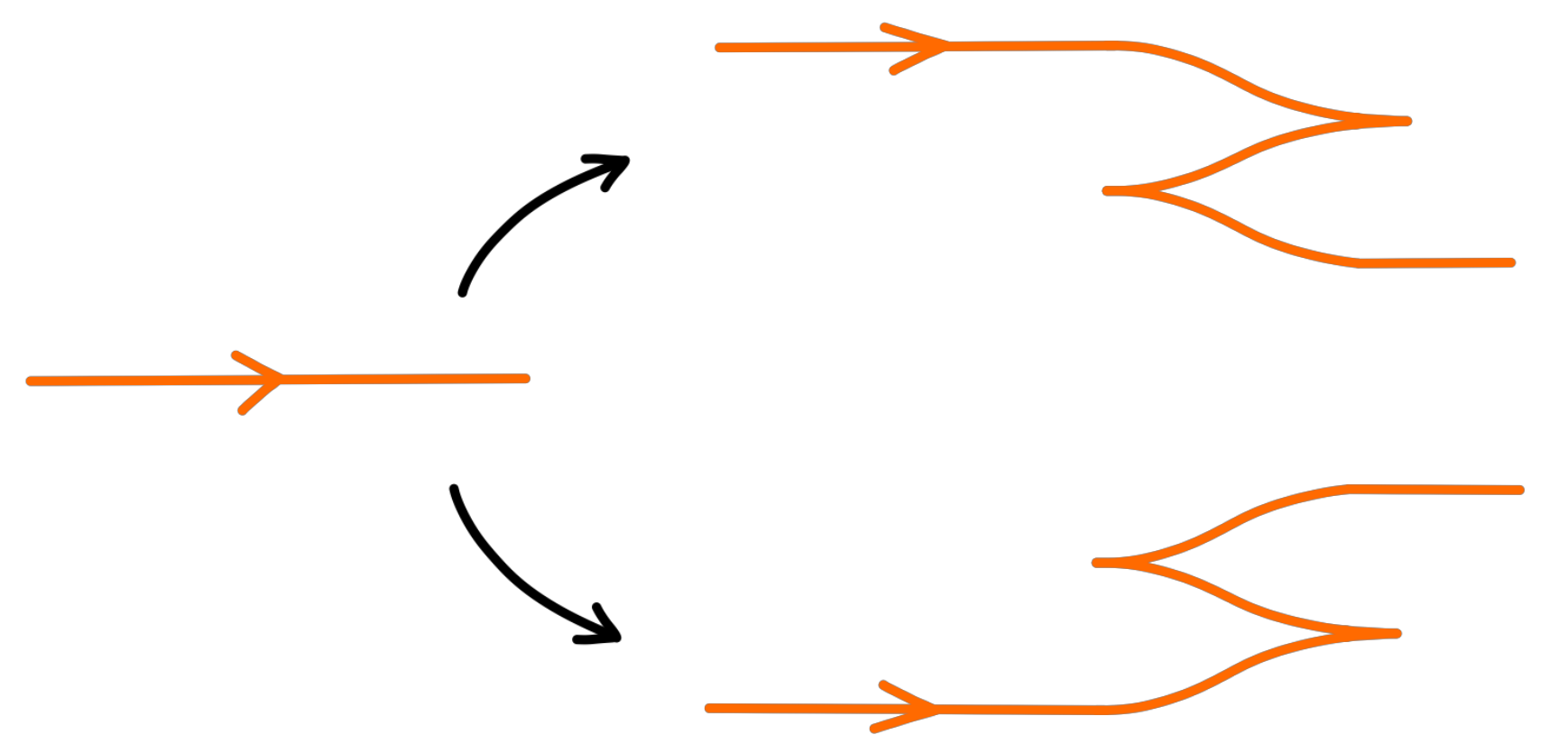}
\caption{Positive and negative stabilizations.}
\label{stabilization}
\end{figure}

Given a Legenrian knot, there is a procedure, called \emph{stabilization}, which has the effect of reducing its Thurston-Bennequin number. After choosing an orientation of the knot, the stabilization can be either positive or negative, according to the effect it has on the rotation number (see figure \ref{stabilization}). 

\noindent Recall that the rotation number of an oriented Legendrian knot in $(S^3,\xi_{st})$ can be computed in the front projection by the formula
\[\rot(K)=\frac{1}{2}(c_D-c_U),\]
where $c_D$ and $c_U$ are respectively the number of down and up cusps. In particular, a positive stabilization increases the rotation number by 1 and a negative stabilization decreases it by 1.

\vspace{0.6cm}

\noindent We put the link of figure \ref{chain} into Legendrian position with respect to the standard tight contact structure of $S^3$, in order to form a linear chain of Legendrian unknots. We do this in such a way that the Thurston-Bennequin number of the $i^{th}$ component is $-a_i+1$. 

The convention that we use here is that  each knot in the front projection is oriented in the counter-clockwise direction.

\begin{figure}[h!]
\centering
\begin{subfigure}[t]{.5\textwidth}
  \centering
  \includegraphics[scale=0.4]{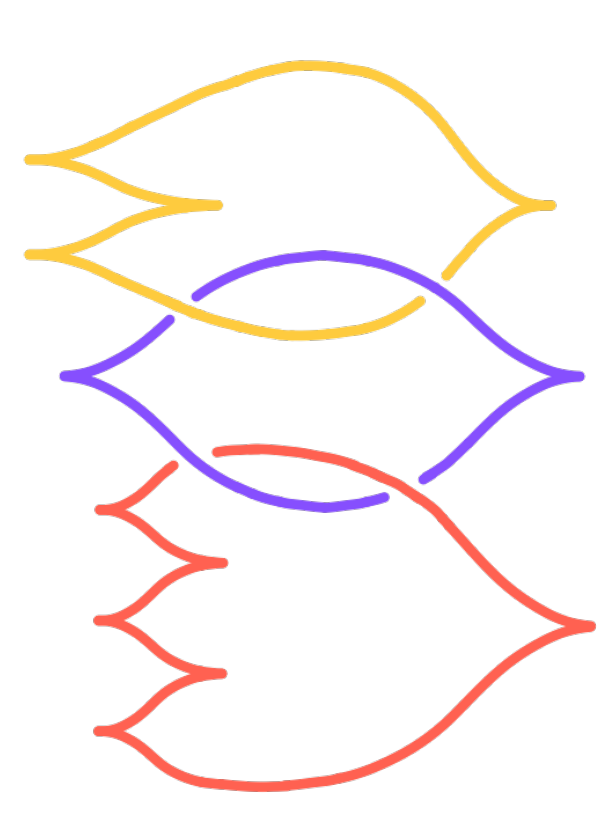}
  \caption{Universally tight structure.}
  \label{egut}
\end{subfigure}%
\begin{subfigure}[t]{.5\textwidth}
  \centering
 \includegraphics[scale=0.4]{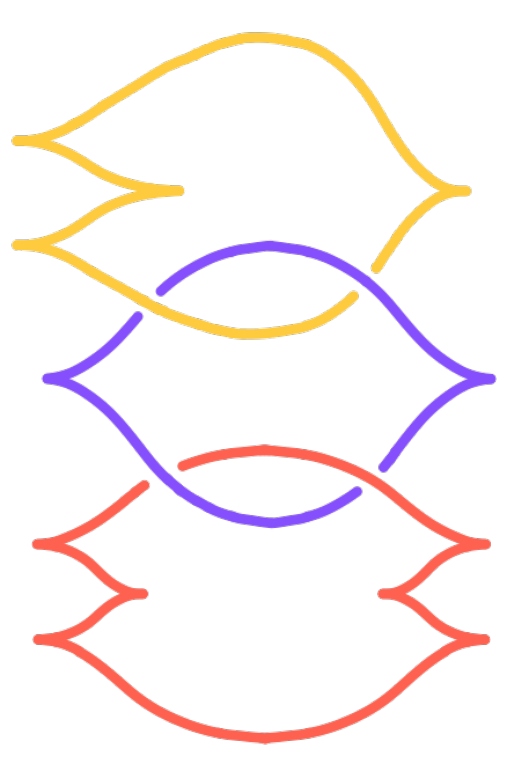}
  \caption{Virtually overtwisted structure.}
  \label{egvot}
\end{subfigure}
\caption{Comparing universally tight and virtually overtwisted contact structures.}
\label{egutvot}
\end{figure}

\noindent The information about the tight contact structure that we get by performing Legendrian (-1)-surgery on this link is encoded by the position of the zig-zags. From the classification of tight contact structures on lens spaces (see \cite{honda}), by looking at the Legendrian link we can tell if the resulting contact structure will be universally tight or virtually overtwisted: if we have only stabilizations of the same type, either all positive or all negative, (i.e. zig-zags on the same side) then the contact structure will be universally tight, otherwise, if we have both positive and negative stabilizations, it will be virtually overtwisted (see figure \ref{egutvot} for an example with $n=3$).

\begin{nrem}
\hypertarget{initialremark}{In} both cases, by attaching $n$ 4-dimensional 2-handles to $B^4$ with framing specified by the Thurston-Bennequin number of each component (decreased by 1), we get a Stein domain whose boundary has the contact structure specified by the Legendrian link. This shows that every tight structure on any lens space is Stein fillable. (We will often write just "filling" to mean "Stein filling").
\end{nrem}

We want to classify the fillings of the virtually overtwisted structures on $L(p,q)$ when $n=2$, i.e. 
\[\frac{p}{q}=[a_1,a_2].\]
The foundational result that makes this possible is the translation from a Legendrian link representing a contact 3-manifold to a presentation via a fibration structure encoded by a compatible open book decomposition (see \cite{giroux}). Loi-Piergallini and independently Akbulut-Ozbagci showed (\cite{loipiergallini}, \cite{akbulut}) that Stein domains can be understood in terms of topological data: a positive allowable Lefschetz fibration over the disk has a Stein structure on its total space, and,  vice versa, any Stein domain can be given such a fibration structure. As a consequence, factorizing the monodromies of \emph{all} the compatible open book decompositions into products of positive Dehn twists is a (theoretical) solution to produce a complete list of Stein fillings for the corresponding contact 3-manifold. In the case of \emph{planar} open book decomposition, the factorization problem is easier thanks to the following theorem of Wendl:

\begin{thm*}[\cite{wendl}]\label{wendl}
If a contact structure $\xi$ on a 3-manifold $Y$ is supported by an open book decomposition with planar page, then every strong symplectic filling of $(Y,\xi)$ is symplectic deformation equivalent to a blow-up of a positive allowable Lefschetx fibration compatible with the given open book.
\end{thm*}

Since all the tight contact structures on lens spaces are planar (see \cite[theorem 3.3]{schonenberger}), we can apply Wendl's theorem to our case.  We will prove theorem \ref{theorem1} by combining techniques coming from mapping class group theory with results by Schönenberger \cite{schonenberger}, Plamenevskaya-Van Horn Morris \cite{plamenVHM}, Kaloti \cite{kaloti} and Menke \cite{menke}.

\section{Proof of the classification theorem}

\begin{defn}
A symplectic filling is called exact if the symplectic form is exact.
\end{defn}

\noindent Stein domains are examples of exact fillings of their boundary. In \cite{menke} it is proved the following:

\begin{thm*} Let $K$ be an oriented Legendrian knot in a contact 3-manifold $(M,\xi)$ and let $(M',\xi')$ be obtained from $(M,\xi)$ by Legendrian surgery on $S_+S_-(K)$, where $S_+$ and $S_-$ are positive and negative stabilizations, respectively. Then every exact filling of $(M',\xi')$ is obtained from an exact filling of $(M,\xi)$ by attaching a symplectic 2-handle along $S_+S_-(K)$.
\end{thm*}

\noindent Then we can derive an immediate corollary (everything is meant up to diffeomorphism).

\begin{cor*} Let $(L,\xi)$ be obtained by Legendrian surgery on the Hopf link.
\begin{itemize} 
\item[a)] Suppose that both components have been stabilized positively and negatively. Then $(L,\xi)$ has a unique Stein filling.
\item[b)] When just one component is positively and negatively stabilized, and the other one has topological framing different from $-4$, then again $(L,\xi)$ has a unique Stein filling.
\item[c)] If only one component is positively and negatively stabilized, and the other one has framing $-4$, then $(L,\xi)$ has two distinct fillings, coming from the two Stein fillings of $(L(4,1),\xi_{st})$.
\end{itemize}
\end{cor*}

\begin{figure}[h!]
  \centering
  \includegraphics[scale=0.3]{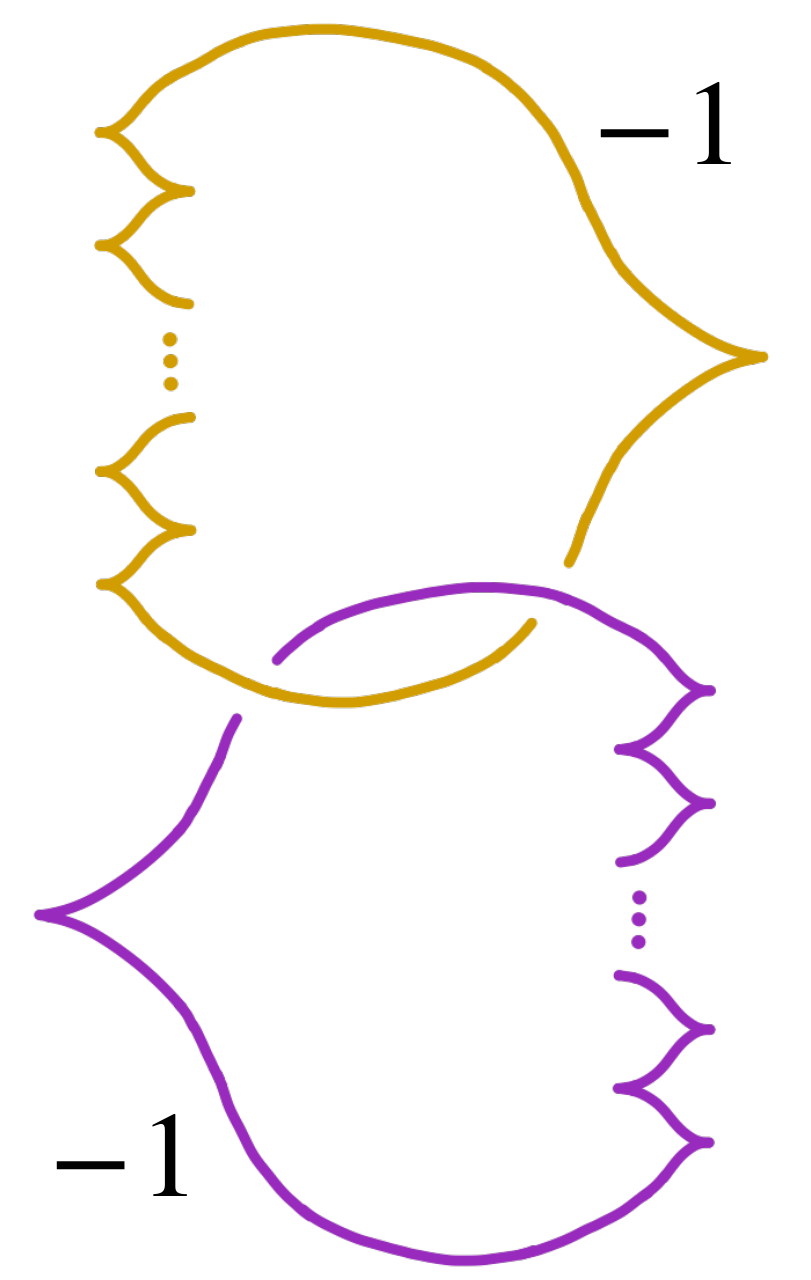}
  \caption{ }
  \label{hopf}
\end{figure}

\noindent The case that does not follow from the theorem of Menke, among the virtually overtwisted structures, is when one component of the link has all the stabillizations on one side and the other component on the opposite side. The rest of the paper is devoted to cover this missing case. We will derive the classification of the fillings of the contact structures on $L(p,q)$ as in figure \ref{hopf} in various steps. 

\noindent Before that, we recall how the stabilization of a Legendrian knot affects the corresponding open book.

\begin{figure}[ht!]
\centering
\includegraphics[scale=0.5]{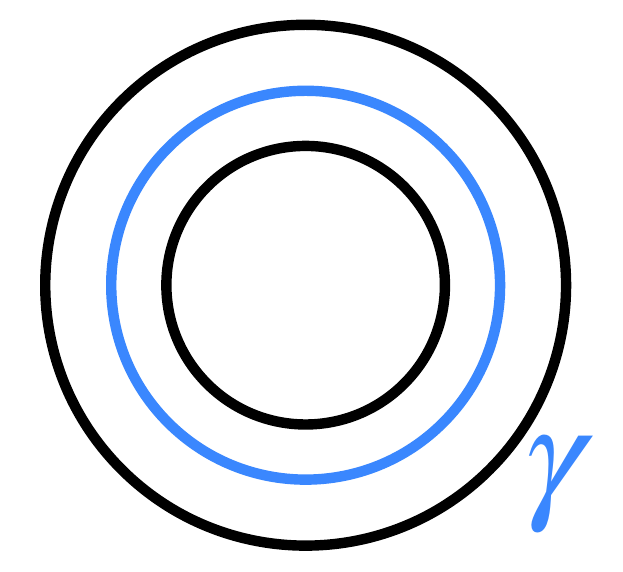}
\caption{Open book decomposition for $(S^3,\xi_{st})$.}
\label{S3page}
\end{figure}

\begin{figure}[h!]
\centering
\begin{subfigure}[t]{.5\textwidth}
  \centering
  \includegraphics[scale=0.5]{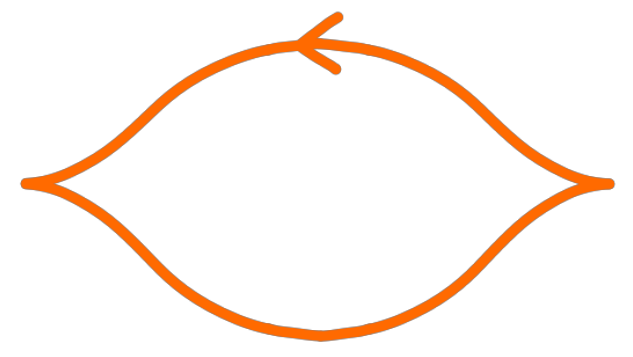}
  \caption{Legendrian unknot with $\tb=-1.$}
  \label{legunknot}
\end{subfigure}%
\begin{subfigure}[t]{.5\textwidth}
  \centering
 \includegraphics[scale=0.5]{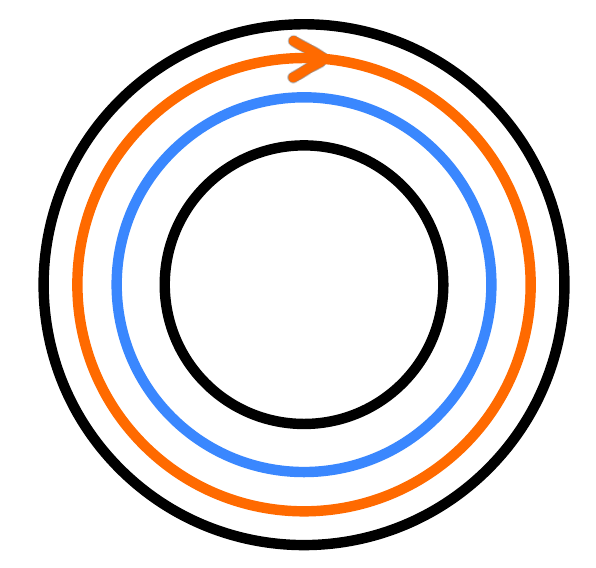}
  \caption{Resulting page.}
  \label{legunknotpage}
\end{subfigure}
\caption{Placing a knot on a page.}
\label{knotannulus}
\end{figure}

We start with the open book decomposition of $(S^3,\xi_{st})$ with annular page and monodromy given by a positive Dehn twist along the core curve $\gamma$, see figure \ref{S3page}. 
%%%% FIGURE OF ANNULUS WITH GAMMA
Now a Legendrian unknot $K$ with $\tb=-1$ in $(S^3,\xi_{st})$ can be placed on a page of the previous decomposition, simply by drawing a parallel copy of $\gamma$ (see \cite{etnyre}). Performing Legendrian surgery on $K$ gives a contact manifold with compatible open book whose page is an annulus and whose monodromy is $\tau_{\gamma}\tau_{\gamma}$ (figure \ref{knotannulus}). 

 %%%% DOUBLE FIGURE OF UNKNOT AND ANNULUS WITH GAMMA TWICE
 
\begin{figure}[h!]
\centering
\begin{subfigure}[t]{.3\textwidth}
  \centering
  \includegraphics[scale=0.5]{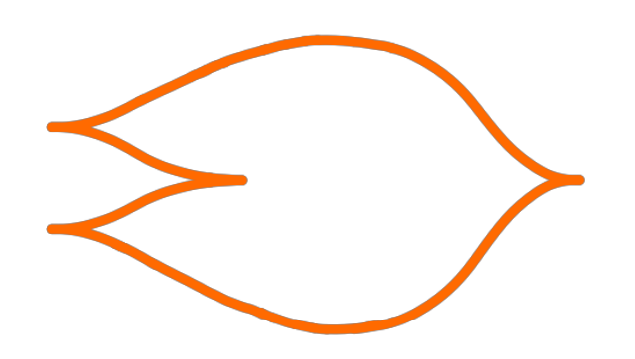}
  \caption{Positive stabilization.}
  \label{posstab1}
\end{subfigure}%
\begin{subfigure}[t]{.3\textwidth}
  \centering
 \includegraphics[scale=0.07]{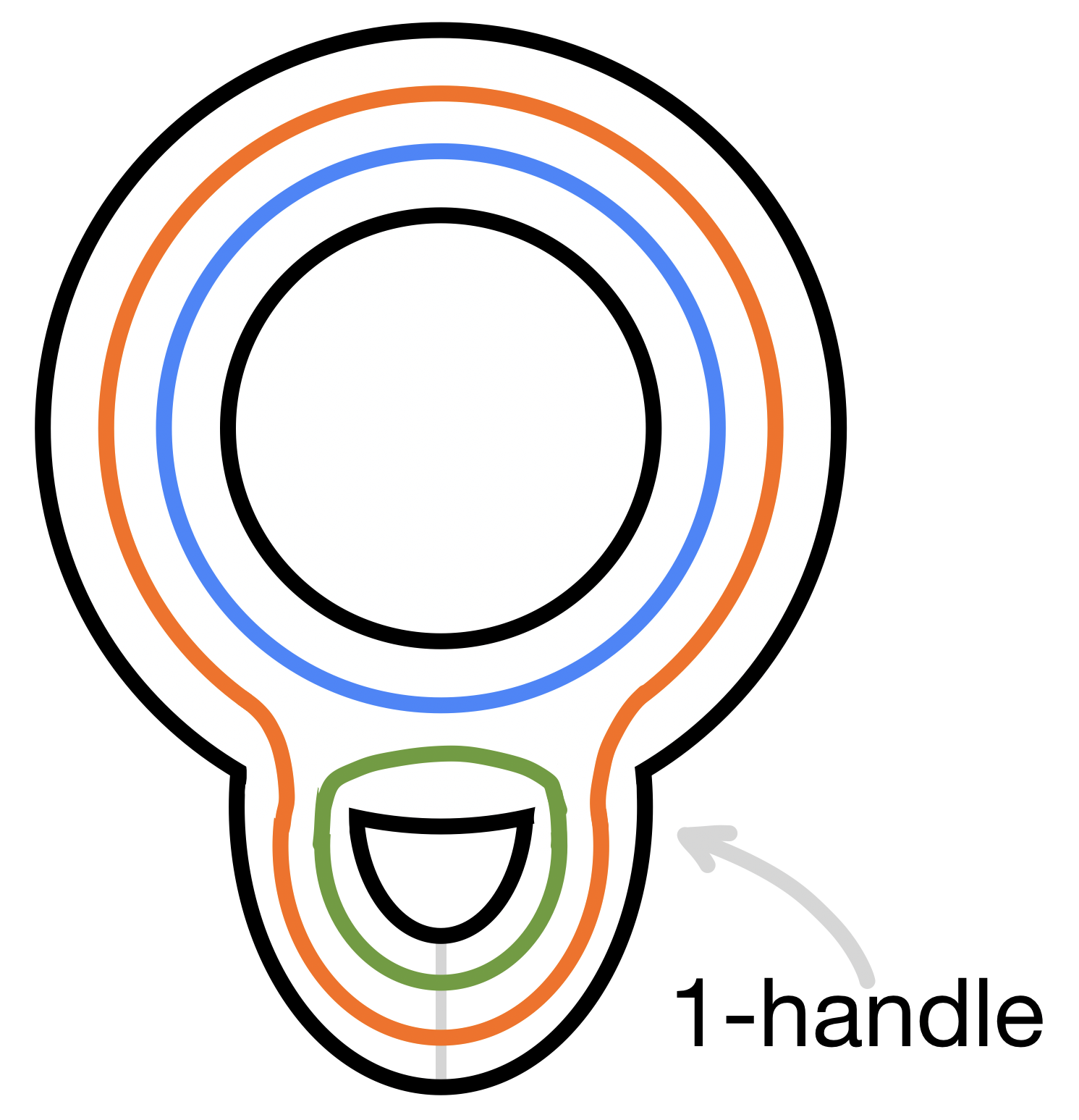}
  \caption{Effect on the open book.}
  \label{posstab2}
\end{subfigure}
\begin{subfigure}[t]{.3\textwidth}
  \centering
 \includegraphics[scale=0.5]{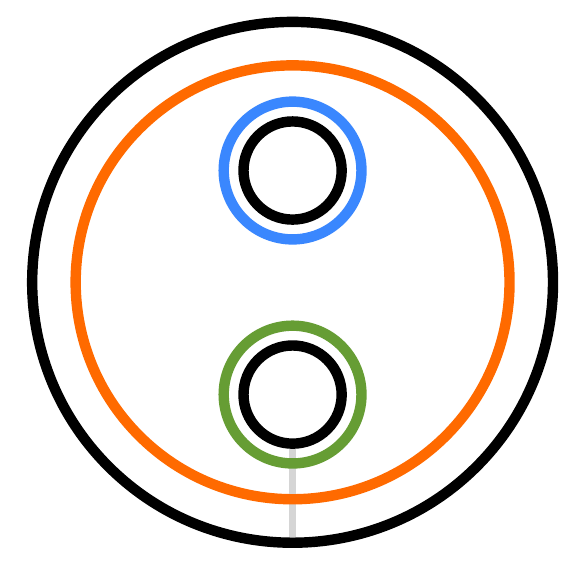}
  \caption{Resulting page.}
  \label{posstab3}
\end{subfigure}
\begin{subfigure}[t]{.3\textwidth}
  \centering
  \includegraphics[scale=0.5]{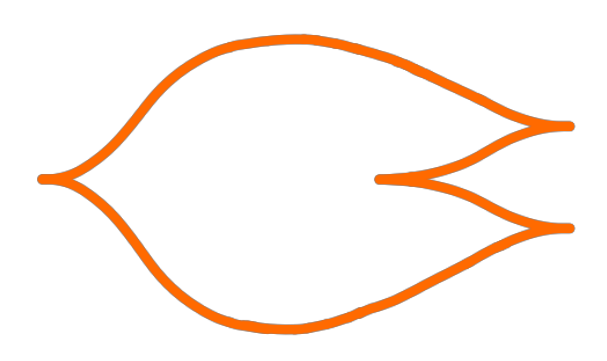}
  \caption{Negative stabilization.}
  \label{negstab1}
\end{subfigure}%
\begin{subfigure}[t]{.3\textwidth}
  \centering
 \includegraphics[scale=0.06]{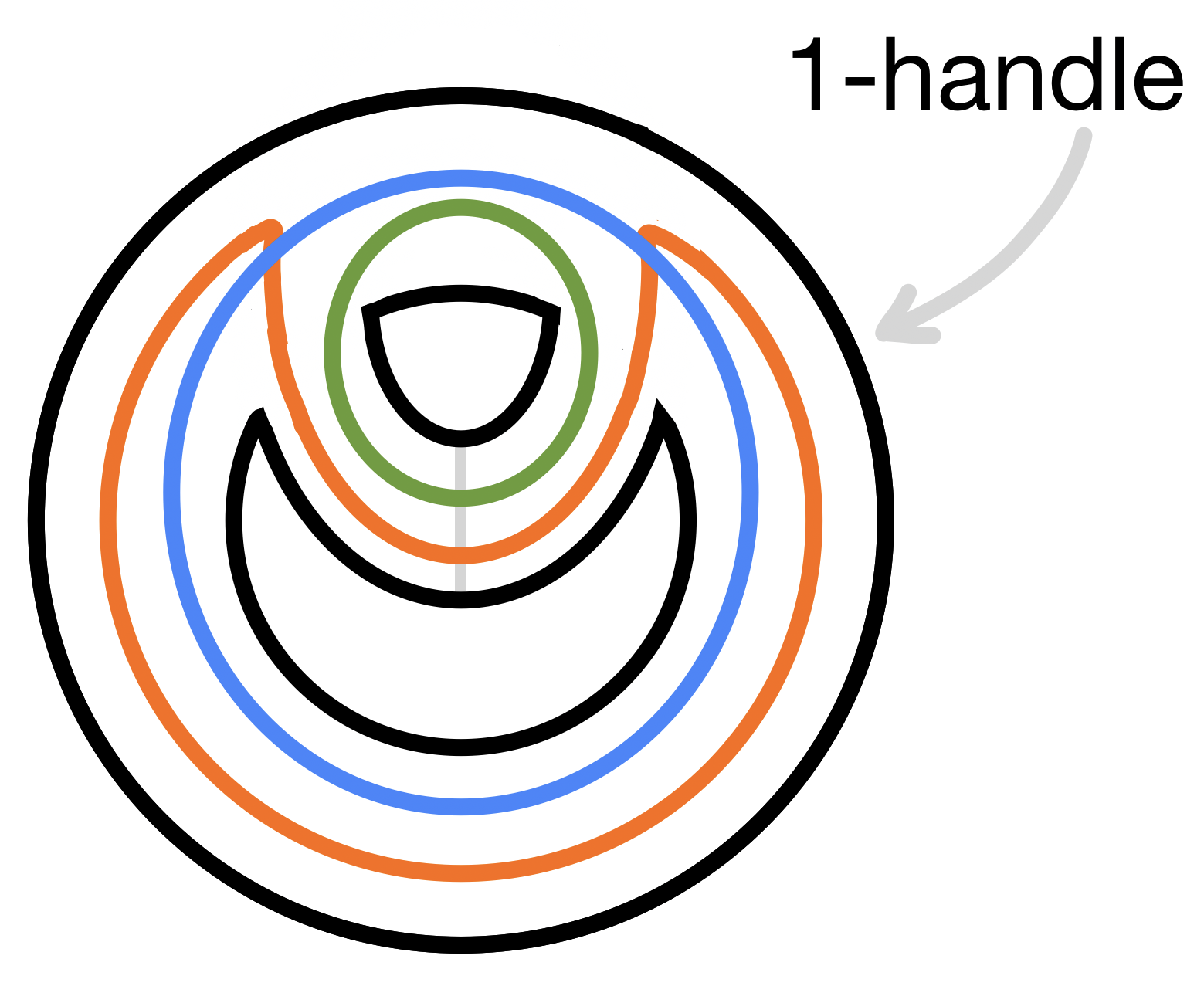}
  \caption{Effect on the open book.}
  \label{negstab2}
\end{subfigure}
\begin{subfigure}[t]{.3\textwidth}
  \centering
 \includegraphics[scale=0.5]{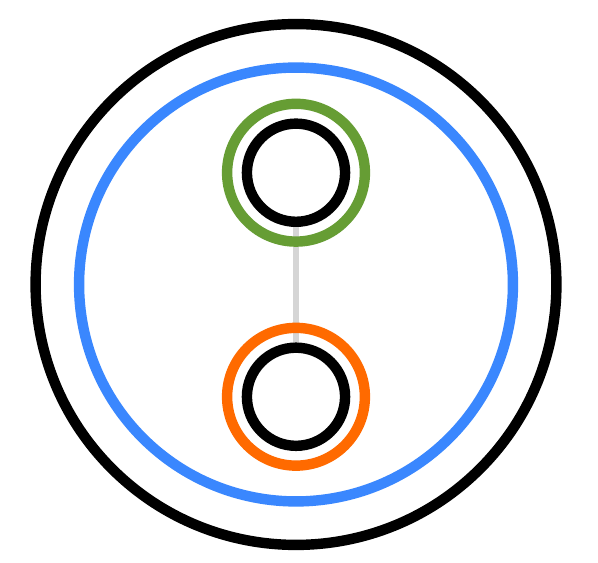}
  \caption{Resulting page.}
  \label{negstab3}
\end{subfigure}
\caption{Placing a stabilized knot on a page.}
\label{posnegstab}
\end{figure}

To reduce the Thurston-Bennequin number of the unknot we add a positive or negative stabilization. Before drawing this Legendrian unknot $K'$ on a page of an open book decomposition of $(S^3,\xi_{st})$ we need to attach a 1-handle to the annulus and modify the monodromy by adding a positive Dehn twist along a curve intersecting the cocore of this new 1-handle once. We always attach the 1-handles on a connected component of the boundary, so that the total number of boundary components increases by one every time. Then we can draw $K'$ on the page by sliding the core curve of the annulus over the 1-handle. According to where we attach the 1-handle, we get a positive or negative stabilization (compare with figure \ref{posnegstab} and with \cite{etnyre}): in order to distinguish between positive and negative, we need to pick an orientation of $K'$ on the page and we orient it in the clockwise direction.

For the rest of the article we will always use this \textbf{orientation convention}, that we stress one more time: Legendrian knots in the front projection are oriented in the counter-clockwise direction and knots on the page of an open book are oriented in the clockwise direction, as in figure \ref{knotannulus}.

 %%%% QUADRUPLE FIGURE OF UNKNOT  WITH POS NEG STAB AND THE 2 OPEN BOOK
 
\subsection*{Step 1: place the link on a planar page}

We will focus on those lens spaces $L=L(p,q)$ with $\frac{p}{q}=[a_1,a_2]$. 
Starting from the link of figure \ref{hopf}, we construct a planar open book decomposition of $(S^3,\xi_{st})$ so that the link itself can be placed on a page inducing the same framing as the contact framing: in order to do so, we first slide one component over the other. This changes the isotopy class of the link, but does not change (see \cite{ding}) the contact type of the 3-manifold obtained by Legendrian (-1)-surgery. 

\begin{figure}[h!]  
  \centering
 \includegraphics[scale=0.5]{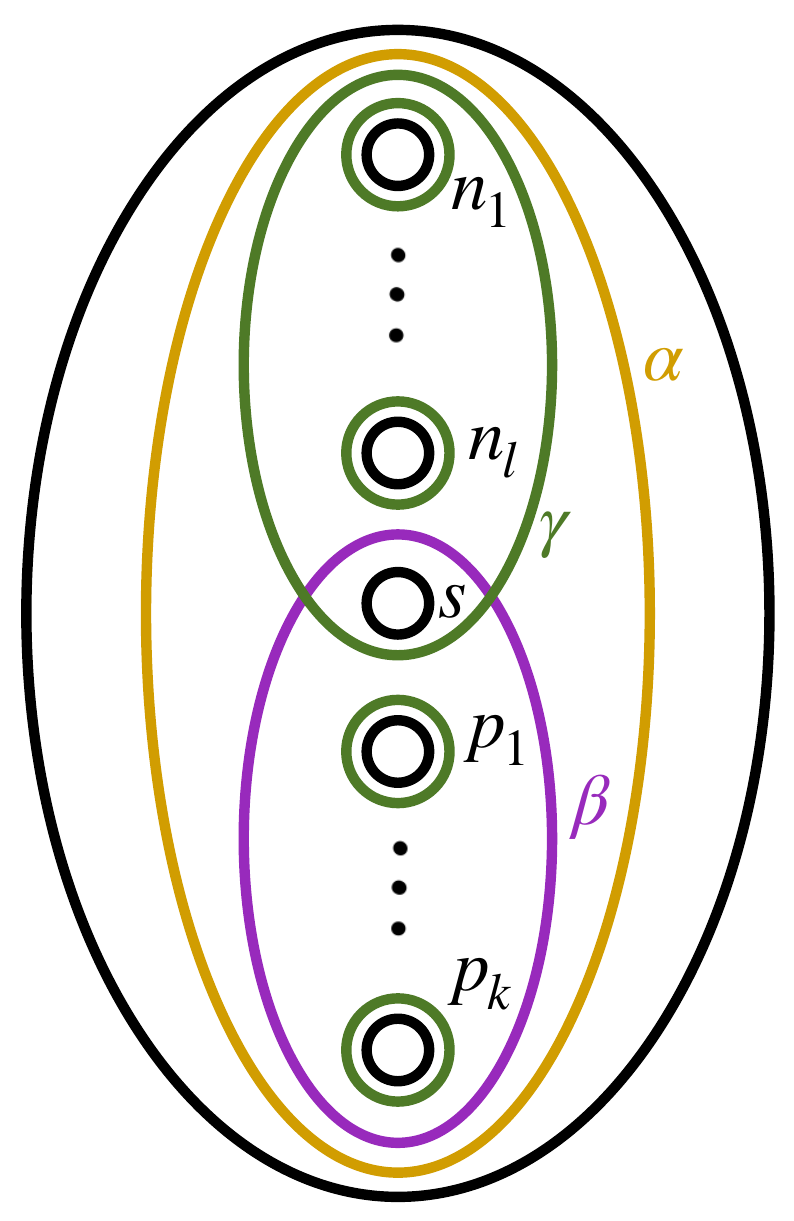}
\caption{Placing the (rolled-up) link on a page.}
\label{page}
\end{figure}

\begin{prop} Figure \ref{page} represents the page (and the monodromy) of an open book decomposition compatible with the contact structure obtained by performing Legendrian surgery on the link of figure \ref{hopf}.
\end{prop}

\begin{prf}
Proceeding as described above, a parallel copy $\alpha$ (yellow) of the core curve $\gamma$ (green) gets positively stabilized as many times as needed, resulting in the addition of these 1-handles to the annulus, on which $\alpha$ itself is slid. We also add the green curves corresponding to monodromy change.

Now an extra curve $\beta$ (purple) parallel to the stabilized $\alpha$ on the page corresponds in $(S^3,\xi_{st})$ to an unknot running parallel to $\alpha$.

By adding copies of the Hopf band we can further stabilize $\beta$ and end up with the desired result, figure \ref{page} (the green curves are the ones coming from the stbilizations plus the starting one descrbing $S^3$). 
\end{prf}

\noindent Call $p_i$ the stabilizing curves corresponding to the positive stabilization of the knot, and $n_i$ the stabilizing curves corresponding to the negative ones.

This is how we can place the link (after sliding) on a planar page of an open book for $(S^3,\xi_{st})$, as described in \cite{schonenberger}. Here is the advantage of such a construction: performing Legendrian surgery on the link $\alpha\cup\beta$ gives the same contact 3-manifold as the one described by the abstract open book decomposition with that page and monodromy given by post-composing the original monodromy with $\t_{\alpha}\t_{\beta}$, that we eventually call $\varphi$. Hence, by Giroux's correspondence we have that the contact type of $(L(p,q),\xi_{vo})$ is encoded in the pair $(\S,\varphi)$ with 
\begin{equation}\label{monod}
\varphi=\t_{\alpha}\t_{\beta}\t_{\gamma}\t_{p_1}\cdots \t_{p_k}\t_{n_1}\cdots \t_{n_l},
\end{equation}
which in turn gives a Stein filling of $(L(p,q),\xi_{vo})$.

The theorem of Wendl that we stated above implies that it is enough to find all the possible factorizations into positive Dehn twists of a \emph{given} planar monodromy in order to get all the Stein fillings of the contact manifold it represents.

\subsection*{Step 2: compute the possible homological configurations}

In their work, Plamenevskaya and Van Horn-Morris \cite{plamenVHM} introduce the multiplicity and joint multiplicity of one or of a pair of holes for a given element in the mapping class group of a planar surface which is written as a product of positive Dehn twists. 

To define these numbers we need the "cap map", which is induced by capping off all but one (respectively two) interior component, while the outer boundary component is never capped. In the first case, we get the mapping class group of the annulus, which is isomorphic to $\Z$, generated by a positive Dehn twist along the core curve. In the second case, we get the mapping class group of a pair of pants, isomorphic to a free abelian group of rank 3, generated by the three positive Dehn twists around each boundary component. By projecting onto the third summand (that comes from the outer boundary component) we get the joint multiplicity around the other two components. We denote by $m(-)$ the multiplicity of a single hole, and by $m(-,-)$ the joint multiplicity of a pair of holes.

\begin{minipage}[c]{.40\textwidth}
\centering
\[\xymatrix{
\Gamma(\Sigma) \ar[r]^-{\mbox{cap}} \ar[rd]_-{m(-)} & \Gamma(\Sigma_{0,2}) \ar[d]^-{\simeq} \\
& \Z
}\]
\end{minipage}%
\hspace{10mm}%
\begin{minipage}[c]{.40\textwidth}
\centering
\[\xymatrix{
\Gamma(\Sigma) \ar[r]^-{\mbox{cap}} \ar@/_/[rdd] _-{m(-,-)} & \Gamma(\Sigma_{0,3}) \ar[d]^-{\simeq} \\
& \Z\oplus\Z\oplus\Z \ar[d]^-{\mbox{pr}_3}\\
& \Z
}\]
\end{minipage}

\begin{figure}[h!]
\centering
\begin{subfigure}[t]{.5\textwidth}
  \centering
  \includegraphics[scale=0.5]{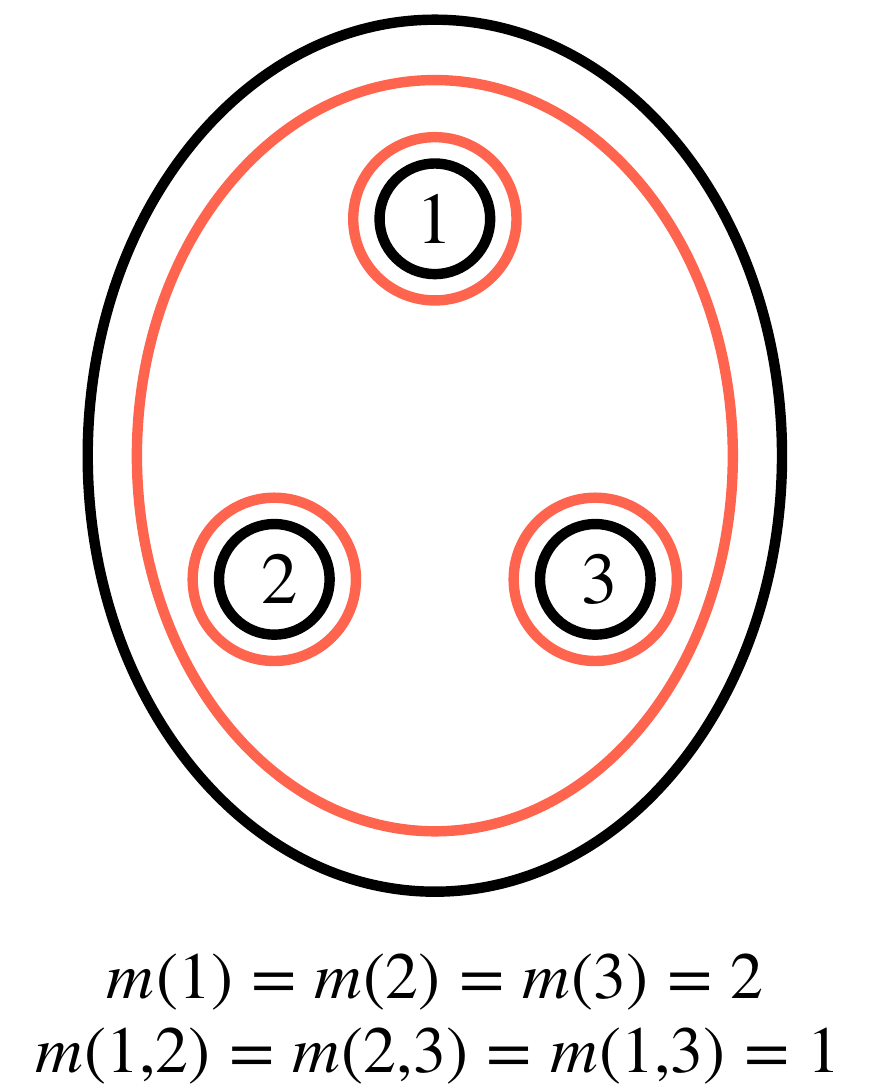}
\end{subfigure}%
\begin{subfigure}[t]{.5\textwidth}
  \centering
 \includegraphics[scale=0.5]{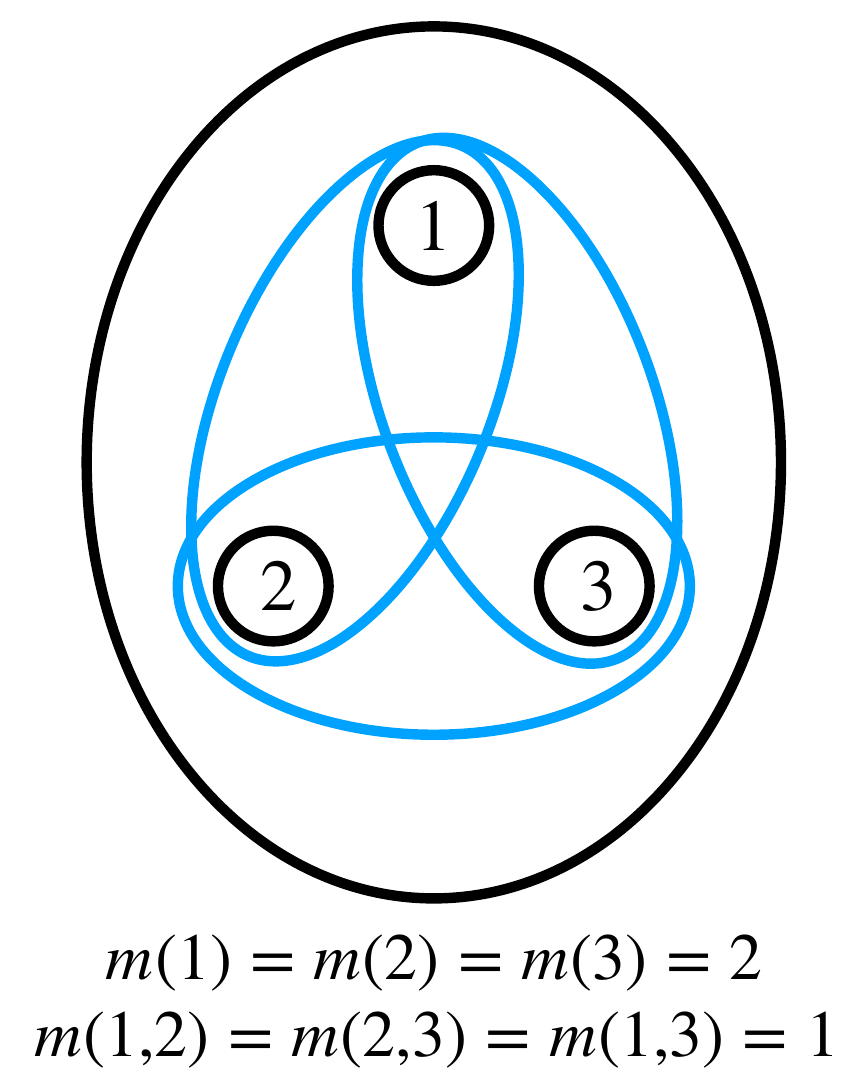}
\end{subfigure}
\caption{Lantern relation and multiplicities.}
\label{lanternmult}
\end{figure}

These multiplicities are independent of the positive factorization we chose for the monodromy.
The reason is that the lantern substitution preserves these numbers (compare with Figure \ref{lanternmult}) and the commutator relation does it too. Since these relations generate all relations in the planar mapping class group (see \cite{margalit}), we can, by introducing negative Dehn twists as well, apply them repeatedly until every curve encloses at most 2 holes.

%%%%% FIGURE OF LANTERN RELATION WITH COMPUTATION OF MULTEPLICITIES

We compute these numbers for the monodromy that we got from Step 1 by looking at Figure \ref{page}. In that figure, the holes called $n_i$'s are the ones corresponding to the negative stabilizations, while the $p_i$'s come from the positive ones and $s$ is the starting hole of the annulus, i.e. the one without a boundary-parallel curve around it. By applying the definition of the cap maps, we directly compute:

\begin{itemize}  % k positive, l negative 
\item $m(n_i)=m(p_i)=3$
\item $m(s)=3$ % starting hole of S^3
\item $m(n_i,n_j)=m(p_i,p_j)=2$
\item $m(n_i,s)=m(p_i,s)=2$
\item $m(n_i,p_j)=1$
\end{itemize}

\noindent Starting from this collection of numbers, we try to reconstruct the homology classes of the curves appearing as support for the positive Dehn twists in the monodromy. We have to distinguish two cases: the first case is when $a_1\neq 4$ and $ a_2\neq 4$ (remember that $[a_1,a_2]$ is the continuous fraction expansion of $p/q$), the second case is when at least one of them is equal to 4. We postpone this second case to \hyperlink{thesentence}{Step 5}. 

\begin{prop}\label{conf}
Assume $a_1\neq 4$ and $ a_2\neq 4$. Then there is a unique homology configuration of curves whose associated multiplicities and double multiplicities are as above, and that gives the same monodromy we started from (Equation \eqref{monod}). In particular, the second Betti number of any Stein filling of the corresponding contact manifold is 2.
\end{prop}

\begin{prf}
We say that a simple closed curve is a \emph{multi-loop} if it encloses at least two holes, in order to distinguish it from a boundary-parallel curve.

Suppose we have another positive factorization of $\varphi$ (see Equation \ref{monod}), and call $\nu_{p_i}$ the number of multi-loops around the hole $p_i$ in this new factorization. 

From the fact that $m(p_i)=m(n_j)=m(s)=3$, we have the upper bound $\nu_{p_i},\nu_{n_j},\nu_s\leq 3$. Moreover, $m(p_i,s)=m(n_j,s)=2$ implies that $\nu_{p_i},\nu_{n_j},\nu_s\geq 2$. Hence we know
\[2\leq\nu_{p_i},\nu_{n_j},\nu_s\leq 3\]
for every $i$ and $j$. We present in details the longest combinatorial part:

\paragraph{Case $a_1,a_2>4$.}
We claim that, in the case $a_1>4$ (which translates into the fact that there are at least three positive holes), there is a multi-loop encircling all the positive holes and $s$. Since $m(p_1,s)=m(p_2,s)=2$ and $m(s)=3$, there must be a curve $a$ encircling $s,p_1,p_2$. If, by contradiction, there exists a hole $p_{\overline{\imath}}$ which is not encircled by $a$, then by the fact that $m(s,p_{\overline{\imath}})=2$, one finds other two curves $b$ and $c$ such that: $b$ encircles $s,p_{\overline{\imath}},p_2$ but not $p_1$, and $c$ encircles $s,p_{\overline{\imath}},p_1$ but not $p_2$. Since $m(s)=3$ and $m(s,n_1)=2$ we must have one of the three curve $a,b,c$ encircling $n_1$ as well, say it is $a$; but then it is impossible to obtain $m(n_1,p_{\overline{\imath}})=1$ without contradicting either $m(s)=3$ or $m(s,n_1)=2$.
This shows that there is a multi-loop $\beta'$ encircling (at least) all the positive holes and $s$. Similarly, there is a multi-loop $\gamma'$ encircling (at least) all the negative holes and $s$.

Now there are two cases:
\begin{itemize}
\item[1)] $\beta'$ and $\gamma'$ coincide, hence the multi-loop $\beta'=\gamma'$ encircles all the hole (i.e. it is parallel to the outer boundary component), and from now on it will be referred to as $\alpha'$;
\item[2)] $\beta'$ and $\gamma'$ are distinct in homology, hence one sees that $\beta'$ cannot encircle negative holes and $\gamma'$ cannot encircle positive holes (this uses the fact that $a_2>4$).
\end{itemize}
We are left to see how we can place the other curves in these two cases in order to get the multiplicities as in previous factorization (Equation \eqref{monod}):
\begin{itemize}
\item[1)] we just forget about $\alpha'$ by lowering all the multiplicities by 1 and then we do again the computation as above. We end up with a curve around the positive holes and $s$ (homologous to $\beta'$), and a curve around the negative holes and $s$ (homologous to $\gamma'$).
\item[2)] We do the same computation as above, by starting with a curve around $s$ and assuming that there is a hole not encircled by it. We get a contradiction with $m(s)=3$. This shows that there must be a curve encircling all the holes, hence parallel to the outer boundary component (i.e. homologous to $\alpha'$).
\end{itemize}
In both cases we get $\nu_s=3$ and $\nu_{p_i}=\nu_{n_j}=2$ for all $i,j$, and we end up with three multi-loops $\alpha'$, $\beta'$ and $\gamma'$, with $\alpha'$ going around all the holes, $\beta'$ around $\{s,n_1,\ldots,n_l\}$, $\gamma'$ around $\{s,p_1,\ldots, p_k\}$. In this way, all the conditions on the joint multiplicities are met, and we just need to add boundary-parallel loops around all the $p_i$'s and $n_j$'s in order to get $m(n_j)=m(p_i)=3$ as required.

\paragraph{Case $a_1=3$ or $a_2=3$.} This case is easier from the combinatorial point of view, and gives the same result. On the other hand, the case $a_1=4$ or $a_2=4$ gives rise to an extra configuration, as discussed later in Step 5.
\end{prf}

This tells us that the homology of any Stein filling for each one of these lens space is fixed. In particular, another factorization of $\varphi$ must be of the form
\[\varphi=\t_{\alpha'}\t_{\beta'}\t_{\gamma'}\t_{p_1}\cdots \t_{p_k}\t_{n_1}\cdots \t_{n_l},\]
where $\alpha',\beta'$ and $\gamma'$ are simple closed curves on $\S$ such that $[\alpha]=[\alpha'], [\beta]=[\beta'], [\gamma]=[\gamma']$ in $H_1(\S;\Z)$. Notice that we do not need to worry about the $p_i$'s and $n_j$'s because the fact that they homologically enclose just one hole implies that they are boundary-parallel, and so their homotopy (and therefore isotopy) class is already determined. Also the homotopy class of $\alpha'$ is determined (since it is boundary-parallel to the outer component) and so $\t_{\alpha}=\t_{\alpha'}$.
Therefore, if the configuration of curves we started from is like the one of Figure \ref{arcs1}, then we already know how to place the boundary-parallel curves appearing in any other factorization (compare with Figure \ref{arcs2}).

In light of this, using the previous factorization of $\varphi$, we see that all the $\t_{p_j}$'s and $\t_{n_i}$'s, together with $\t_{\alpha}$ and $\t_{\alpha'}$, cancel out, leaving us with:
\[\t_{\beta}\t_{\gamma}=\t_{\beta'}\t_{\gamma'}.\]
This relation holds in $\Gamma_{0,k+l+2}=\Gamma(\S_{0,k+l+2})$, and we are asking ourselves if there can be a pair of curves $\{\beta',\gamma'\}$ on $\S$ with the homological condition that $[\beta']=[\beta]$ and $[\gamma']=[\gamma]$ inside $H_1(\S_{0,k+l+2};\Z)$, and such that the product of the corresponding Dehn twists is isotopic to $\t_{\beta}\t_{\gamma}$. We will reduce the problem from $\S_{0,k+l+2}$ to $\S_{0,4}$.

\subsection*{Step 3: reduce the number of boundary components} 

\begin{defn} A diffeomorphism $f:\Sigma\to \Sigma$ which restricts to the identity on $\partial\Sigma$ is \emph{right-veering} if $f(\eta)$ is to the right of $\eta$ at its starting point, for every properly embedded oriented arc $\eta\subseteq\Sigma$, isotoped to minimize $\eta\cap f(\eta)$.
\end{defn}

\begin{figure}[ht!]
\centering
\begin{subfigure}[t]{.5\textwidth}
  \centering
  \includegraphics[scale=0.5]{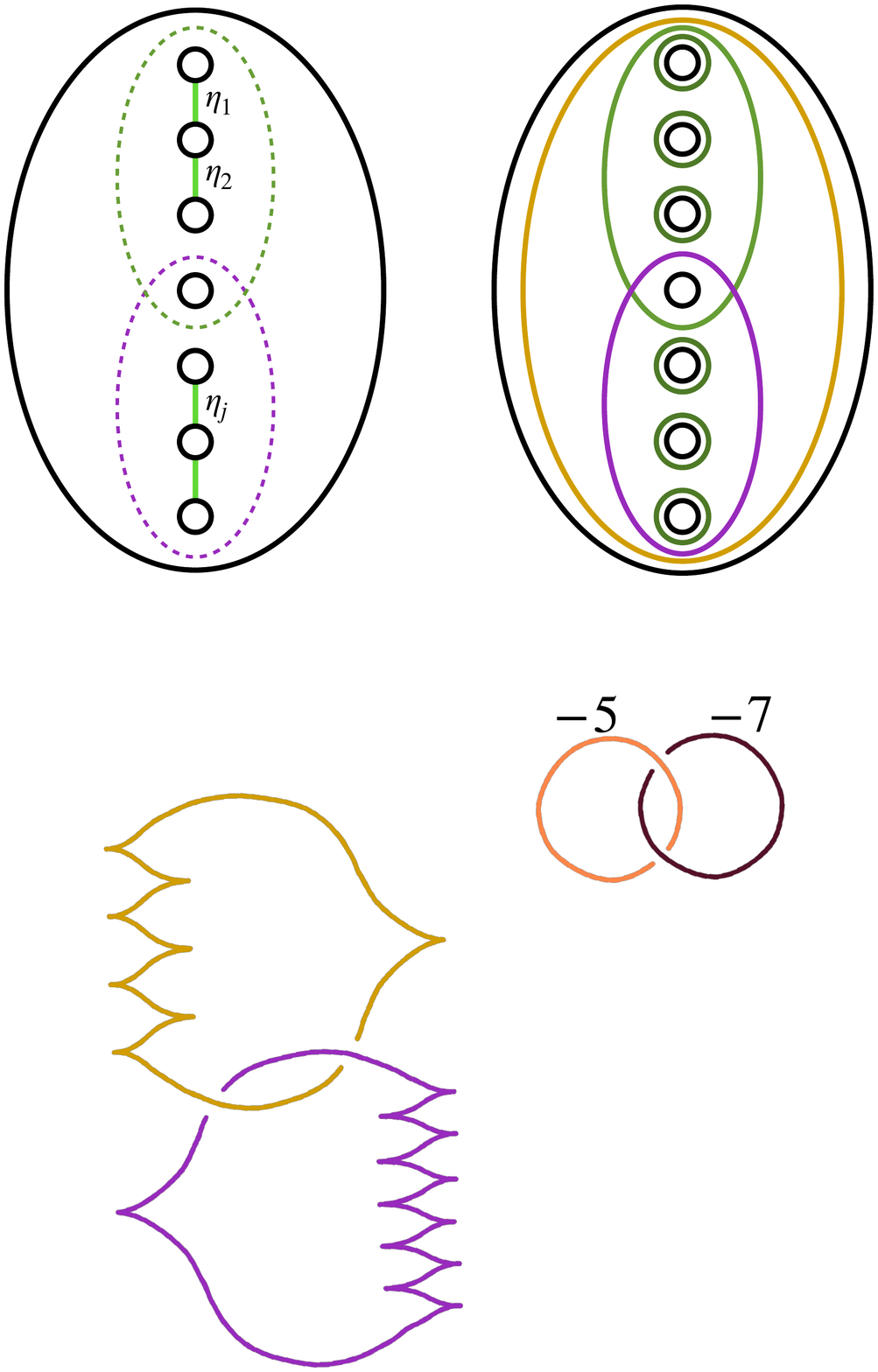}
\caption{Original configuration of curves.}
\label{arcs1}
\end{subfigure}%
\begin{subfigure}[t]{.5\textwidth}
\centering
 \includegraphics[scale=0.5]{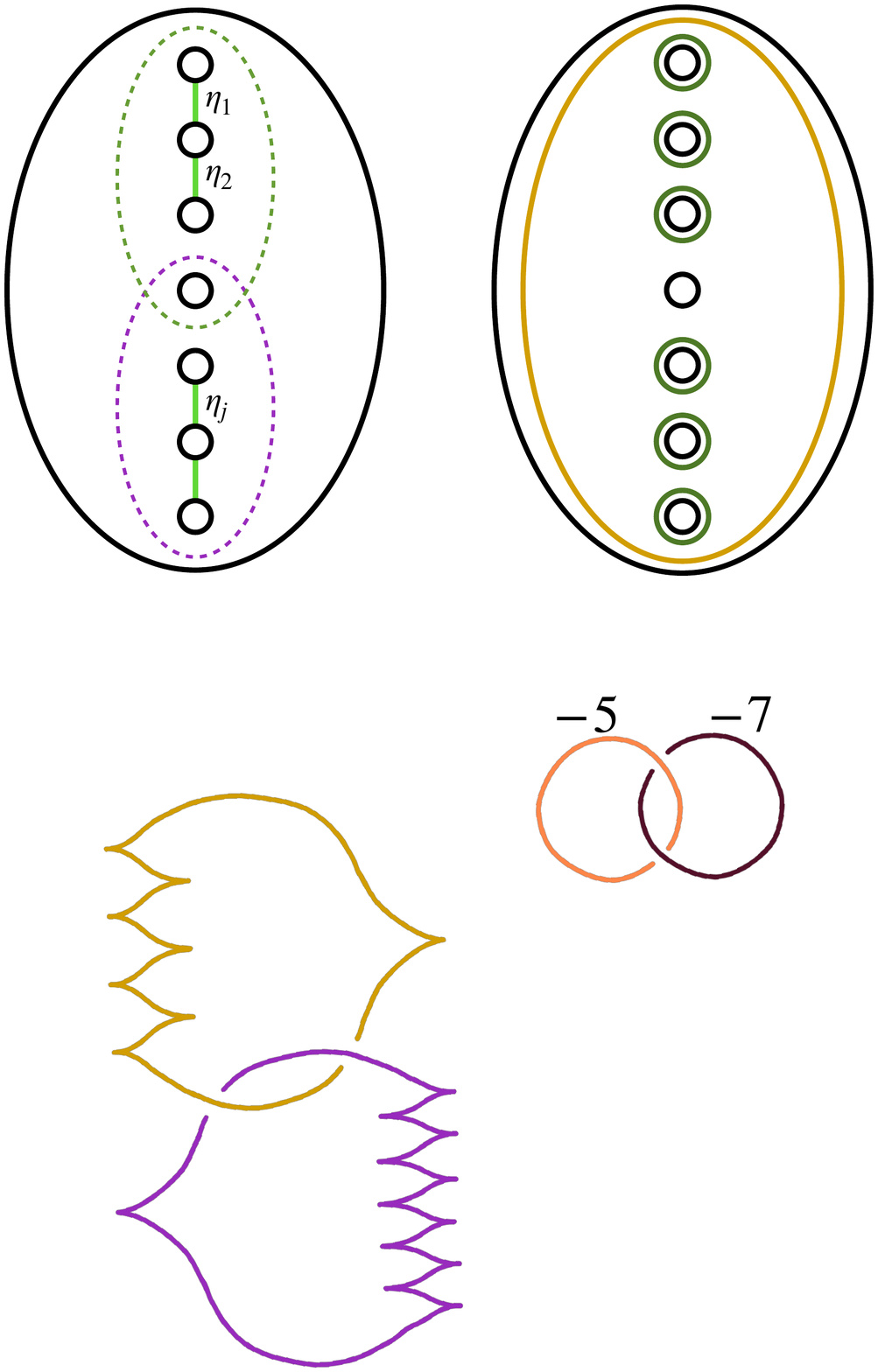}
\caption{Boundary parallel curves.}
\label{arcs2}
\end{subfigure}
\begin{subfigure}[t]{.5\textwidth}
  \centering
  \includegraphics[scale=0.5]{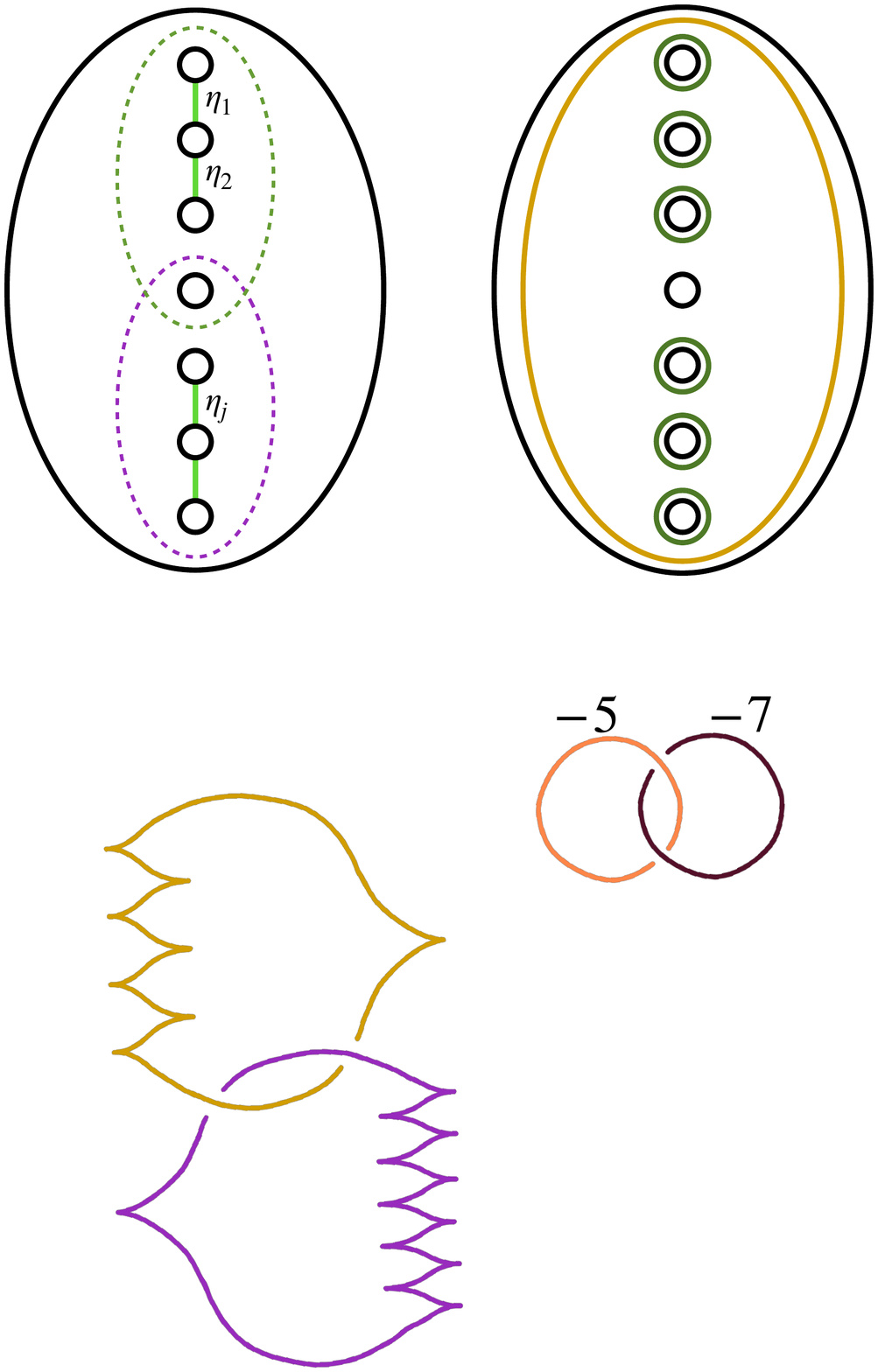}
\caption{Cutting the surface along arcs.}
\label{arcs3}
\end{subfigure}%
\begin{subfigure}[t]{.5\textwidth}
  \centering
 \includegraphics[scale=0.6]{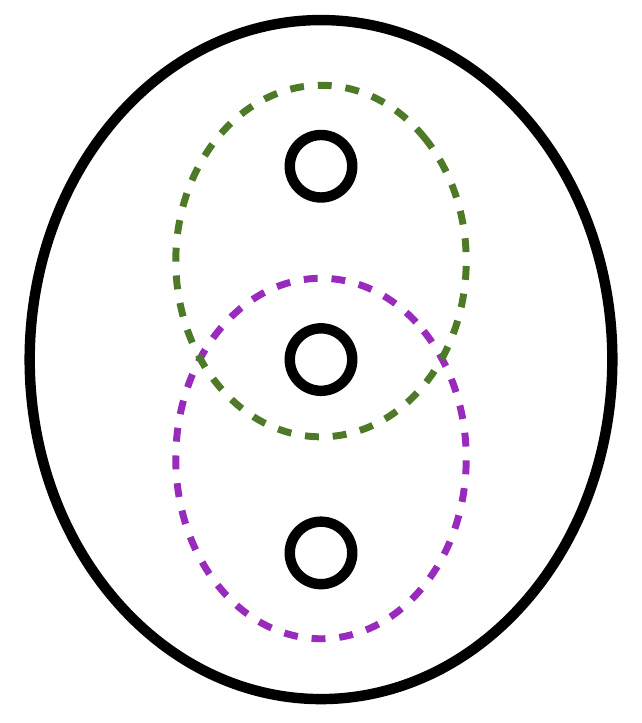}
\caption{Resulting surface and curves.}
\label{arcs4}
\end{subfigure}

\caption{Configuration of curves.}
\end{figure}

\noindent In \cite{rightveering} it is proved that any positive Dehn twist is right-veering, and that the composition of right-veering homeomorphisms is still right-veering. 
Consider the green arcs $\eta_j$ drawn in figure \ref{arcs3}. The dashed curves are drawn like that just to indicate their homology class, while the isotopy classes are still unknown.

The arcs are disjoint from $\beta\cup\gamma$, therefore $\t_{\beta}\t_{\gamma}(\eta_j)=\eta_j$. If any of the curves $\{\beta',\gamma'\}$ (respectively purple and green) crosses one of the $\eta_j$'s, say $\gamma'$, then we would need $\t_{\beta'}$ to move the arc back to the initial position, since 
\[\eta_j=\t_{\beta}\t_{\gamma}(\eta_j)=\t_{\beta'}\t_{\gamma'}(\eta_j).\]
But this is impossible because $\t_{\beta'}$ would move the arc further right, because of right-veering property. Therefore all the arcs are disjoint from $\beta'\cup\gamma'$ as well. We can then cut along these arcs and look for a configuration of $\{\beta',\gamma'\}$ on the resulting surface, which still has genus zero but now has just 4 boundary components (compare with figure \ref{arcs4}).

Notice furthermore that $\alpha$ and $\alpha'$ (pictured in yellow) are parallel to the boundary, so their isotopy class is determined by the homology class, and we can simplify the corresponding Dehn twists on the two side of the equation, leaving us with:
\[\t_{\beta}\t_{\gamma}=\t_{\beta'}\t_{\gamma'} \in \Gamma_{0,4}.\]

\subsection*{Step 4: identify the possible geometric configurations} 

To conclude the argument for the uniqueness of the filling we need a result which is proved by Kaloti:

\begin{lem*}[\cite{kaloti}] Suppose there are two simple closed curves $\beta',\gamma'$ on $\S_{0,4}$ with $[\beta']=[\beta]\in H_1(\Sigma_{0,4};\Z)$, $[\gamma']=[\gamma]\in H_1(\Sigma_{0,4};\Z)$ and such that $\t_{\beta}\t_{\gamma}=\t_{\beta'}\t_{\gamma'} \in \Gamma_{0,4}$. 
Then there exists a diffeomorphism $\Gamma_{0,4}$ taking $\beta\mapsto \beta'$ and $\gamma\mapsto\gamma'$. 
\end{lem*}

\noindent The consequence of this lemma is that the two curves $\beta'$ and $\gamma'$ are, up to diffeomorphism, the same as $\beta$ and $\gamma$, and therefore the filling that the pair $\{\beta',\gamma'\}$ describes, together with the boundary parallel curves, is the one we constructed in the initial \hyperlink{initialremark}{remark}.

\subsection*{Step 5: the extra filling when $a_i=4$ for $i\in\{1,2\}$}\label{extra}

\hypertarget{thesentence}{When} $a_2$ (or $a_1$) is equal to 4, then, as expected, there is another homology configuration which is coherent with the single and joint multiplicities computed above. It is given by applying the lantern relation to the original configuration. This has the effect of reducing by one the total number of curves appearing in the factorization (hence $b_2$ of the corresponding filling is 1 and not 2). 

\begin{prop} 
If $a_2= 4$, then there are two possible homology configurations of curves respecting the single and double multiplicities computed above. 
\end{prop}

\begin{prf} 

If we go through the computation of homological configurations of curves as we did in Proposition \ref{conf}, we find this time two of them, due to the possibility of performing a lantern substitution once (notice in fact that the case when $a_1=a_2=4$ allows again just two configurations, and not three, because applying the first substitution changes the configuration preventing the second-one from being possible).

One configuration has been already described element-wise in Proposition \ref{conf} and corresponds to a Stein filling with $b_2=2$ (unique up to diffeomorphism). 

The other one is homologous (curve by curve) to the configuration we get after applying the lantern substitution to Figure \ref{page} with $l=2$ (i.e. $a_2=4$); the proof of its homology uniqueness is derived as in the proof of Proposition \ref{conf}, and it is omitted here.
\end{prf}

Then we proceed as above: all the boundary-parallel curves are placed and ignored, since, again, their homology classes determine their isotopy classes. Therefore, we can cut along appropriate arcs and reduce the number of holes appearing in the factorization, which, in turn, is the same thing as starting with a Legendrian knot whose Thurston-Bennequin number is smaller. So we can focus on the minimal possible example (after having cut along the maximal system of arcs), which has $[a_1,a_2]=[3,4]$, producing $L(11,4)$, by the computation $3-\frac{1}{4}=\frac{11}{4}$, with the contact structure specified by Figure \ref{hopf(3,4)} whose compatible open book decomposition has page as in Figure \ref{pagehopf(3,4)}. Here we immediately see that we can (uniquely) apply the lantern relation on the set of four curves given by the yellow one, the red one and the two boundary-parallel green curves (compare with Figure \ref{pagelantern}). After the substitution we get a Stein filling with $\chi=2$.

\begin{figure}[ht!]
\centering
\begin{subfigure}[t]{.45\textwidth}
  \centering
 \includegraphics[scale=0.4]{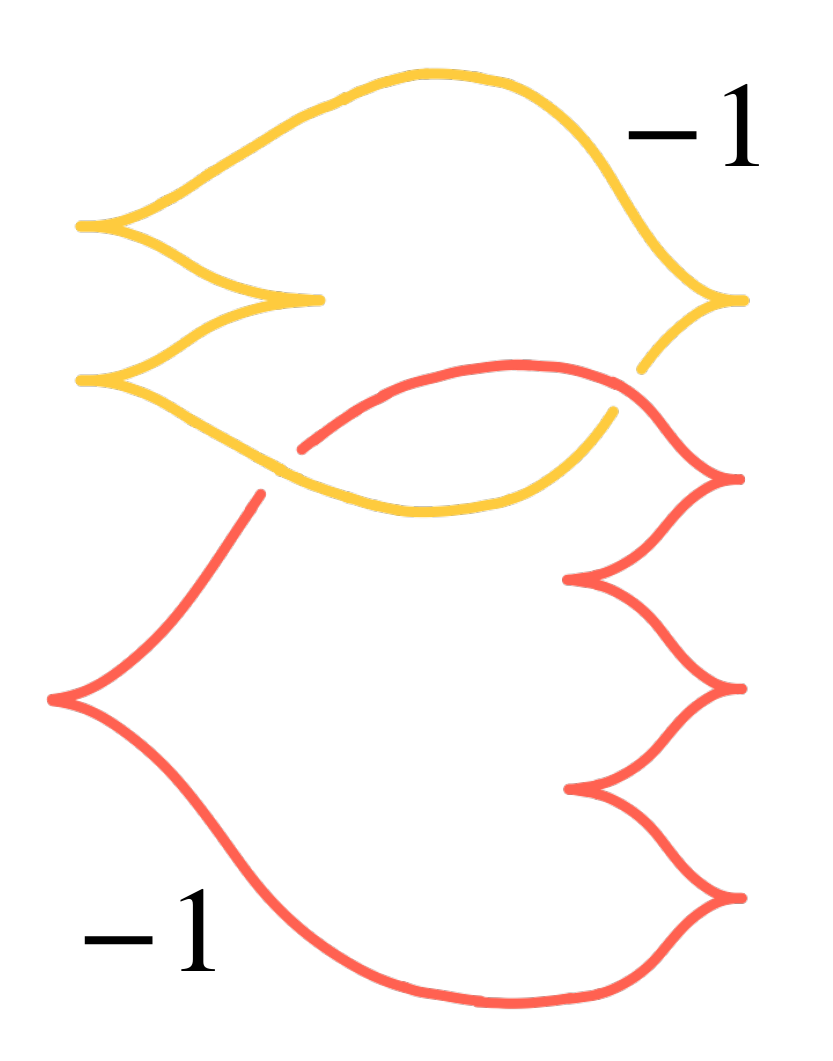}
\caption{ }
  \label{hopf(3,4)}
\end{subfigure}
\begin{subfigure}[t]{.45\textwidth}
  \centering
 \includegraphics[scale=0.3]{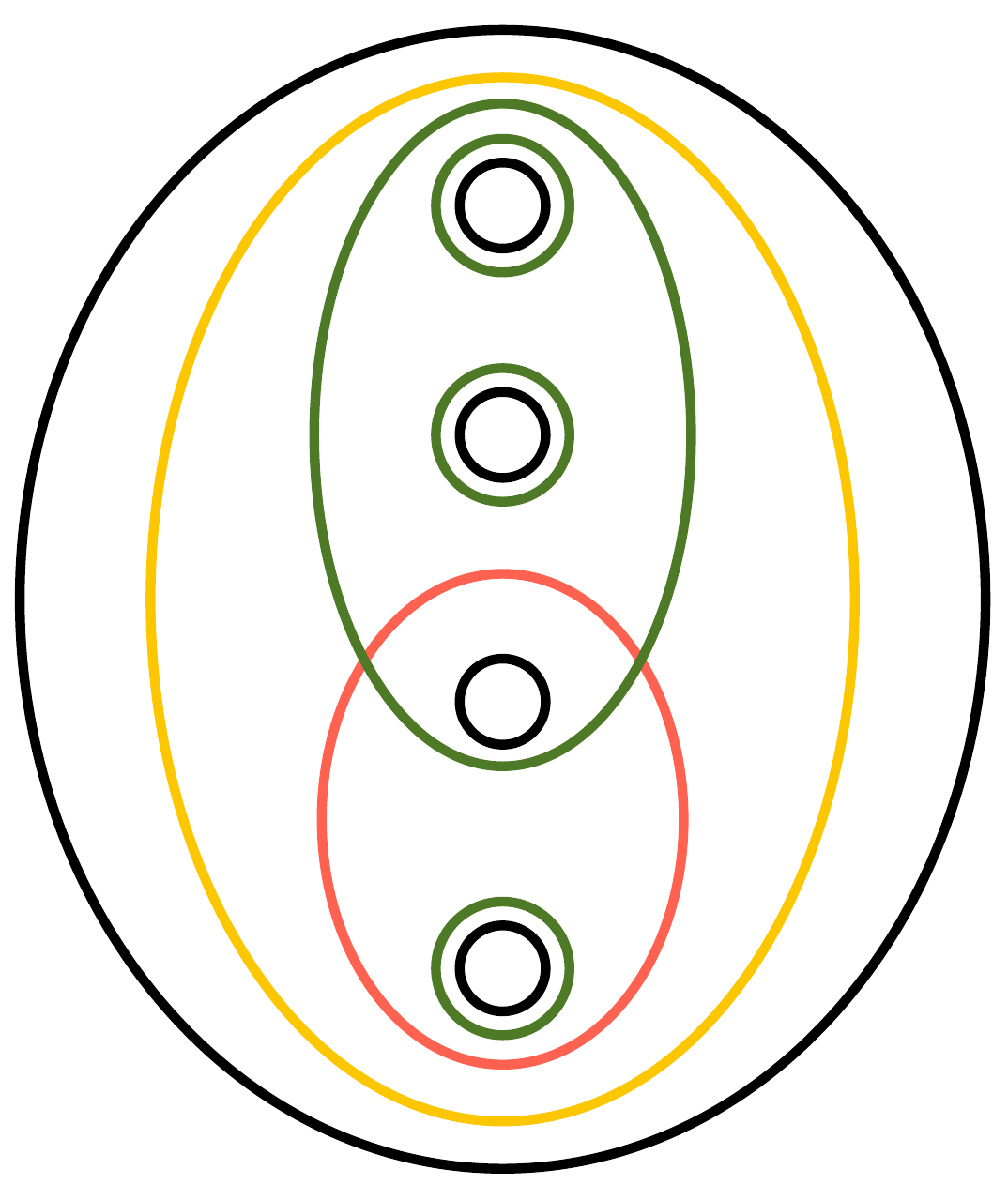}
  \caption{ }
  \label{pagehopf(3,4)}
\end{subfigure}
\caption{Tight structure on $L(11,4)$ with two fillings.}
\label{contactxi}
\end{figure}
 
\begin{prop} Let $X$ be a Stein filling of the contact 3-manifold described by Figure \ref{hopf(3,4)}, with $\chi(X)=2$. Then the homeomorphism type of $X$ is unique.
\end{prop}

\begin{prf} In order to derive our statement we need three facts:

\begin{itemize}
\item[1)] since we have determined the homology configuration of curves appearing in the factorization of the monodromy, we can use Figure \ref{pagelantern} to compute $H_1(X)$: the holes of the surface correspond to the 1-handles of $X$, and the curves themselves are the attaching circles of the 2-handles. It is immediate to see that $H_1(X)=0$, hence $\pi_1(X)$ is perfect. 
But $\pi_1(X)$ is a quotient of $\Z/11\Z\simeq \pi_1(L(11,4))$, see \cite[page 216]{ozbagci}, hence it is abelian. We conclude that $X$ is simply connected.

\item[2)] The intersection form of $X$ is characterized by the homological configuration of curves in the open book decomposition, hence it is uniquely determined (and it is isomorphic to $[-11]$).

\item[3)] The fundamental group of the boundary of $X$ is $\pi_1(L(11,4))\simeq \Z/11\Z$.
\end{itemize}

\noindent Then \cite[Proposition 0.6]{boyer} applies and tells that $X$ is unique up to homeomorphism, as claimed. The reason why we cannot describe all the (potential) diffeomorphism types is because, after applying the lantern substitution, we cannot solve explicitly the geometric configuration problem on $\Sigma_{0,5}$ with the curves we got (see Figure \ref{pagelantern}), and, moreover, there is no arc which is disjoint from these curves and on which to cut open in order to reduce the number of boundary components.
\end{prf}

\section{Final remark}
\begin{figure}[ht!]
\centering
\begin{subfigure}[t]{.45\textwidth}
  \centering
 \includegraphics[scale=0.35]{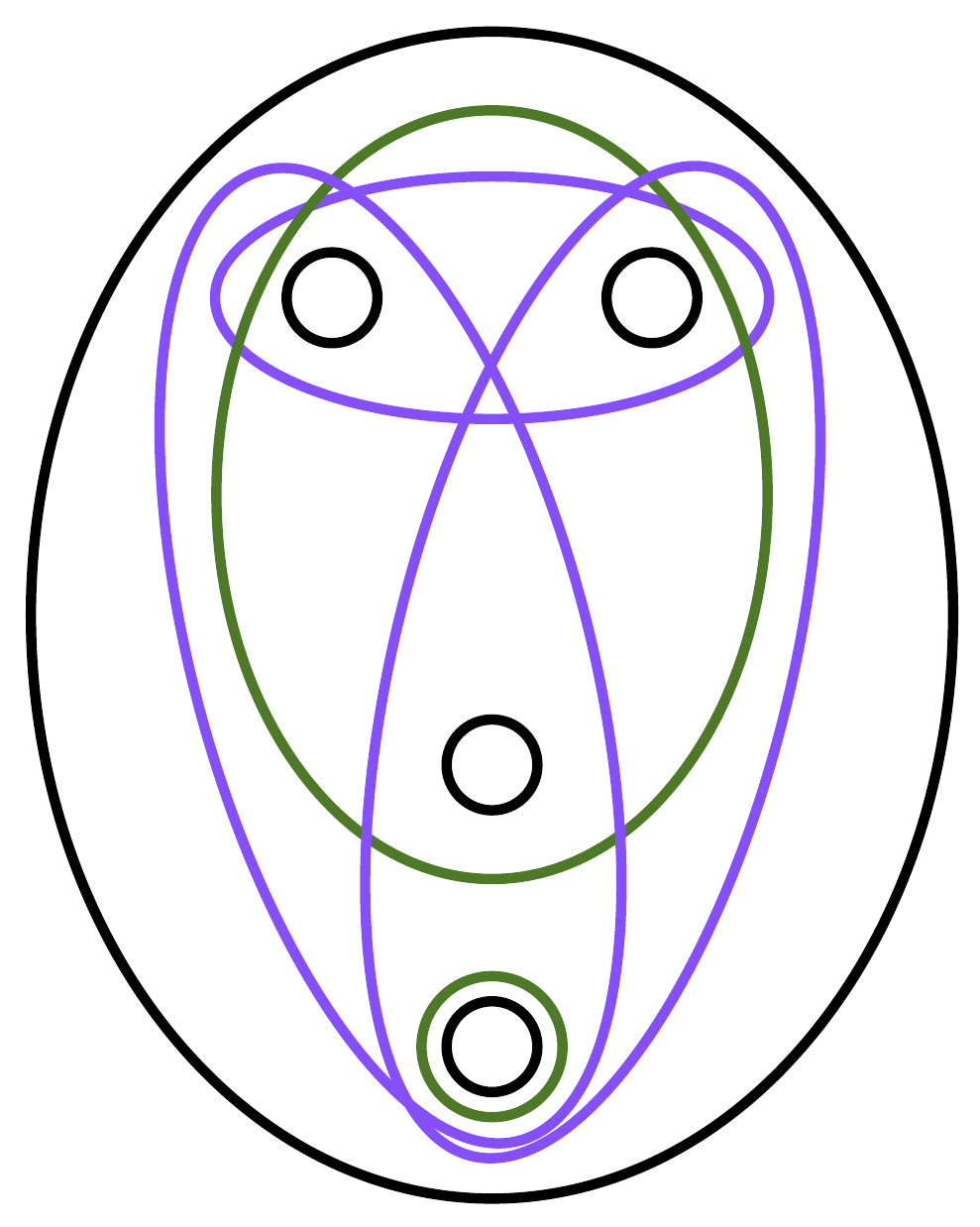}
\caption{The configuration of curves after lantern substitution.}
  \label{pagelantern}
\end{subfigure}
\begin{subfigure}[t]{.45\textwidth}
  \centering
 \includegraphics[scale=0.40]{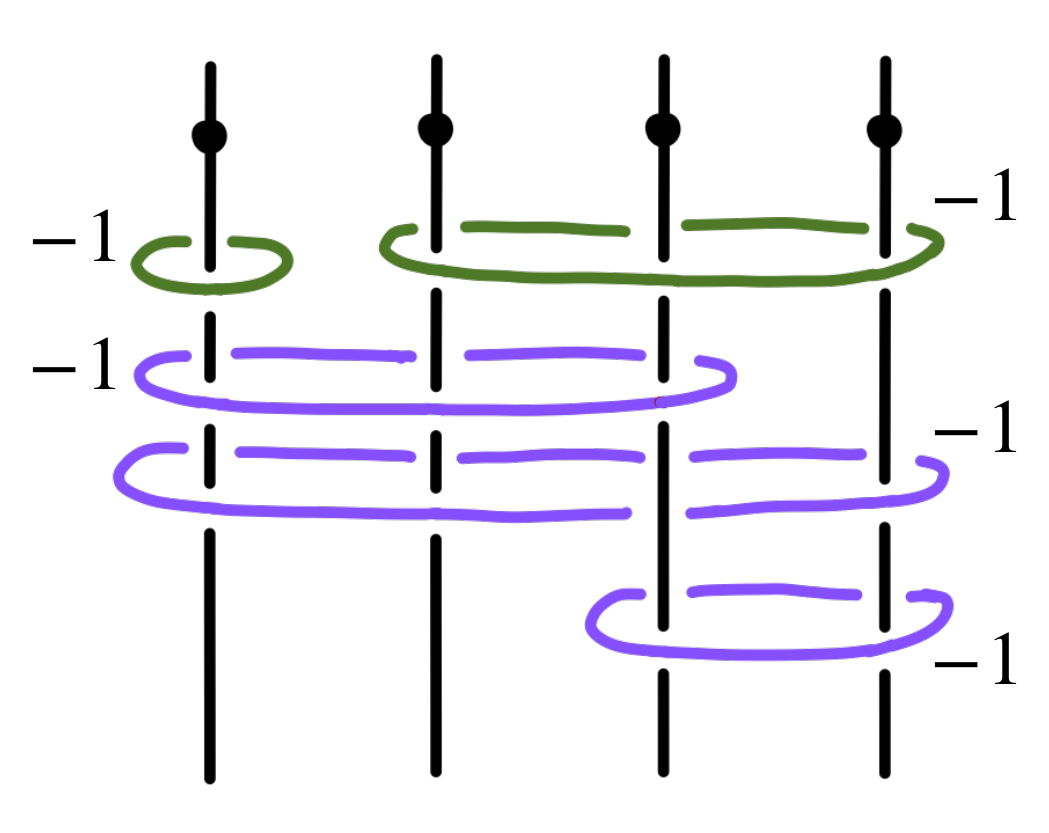}
  \caption{The corresponding 4-manifold X.}
  \label{steinlantern}
\end{subfigure}
\caption{Description of the filling with $\chi=2$.}
\end{figure}

Starting from the explicit configuration of curves corresponding to $L(11,4)$ with the contact structure of Figure \ref{hopf(3,4)}, we can apply the lanter substitution (see Figure \ref{pagelantern}) and then draw a Kirby diagram of the corresponding Stein domain $(X,J)$, see Figure \ref{steinlantern}. 

By performing handle calculus we get a new diagram of $X$ which is simpler in the following sense: this smooth handle decomposition of $X$ consists of a 0-handle and a single 2-handle, attached along the torus knot of type $(-5,2)$, pictured in Figure \ref{t(5,2)}, with framing $-11$.

\begin{figure}[ht!]
\centering
\includegraphics[scale=0.5]{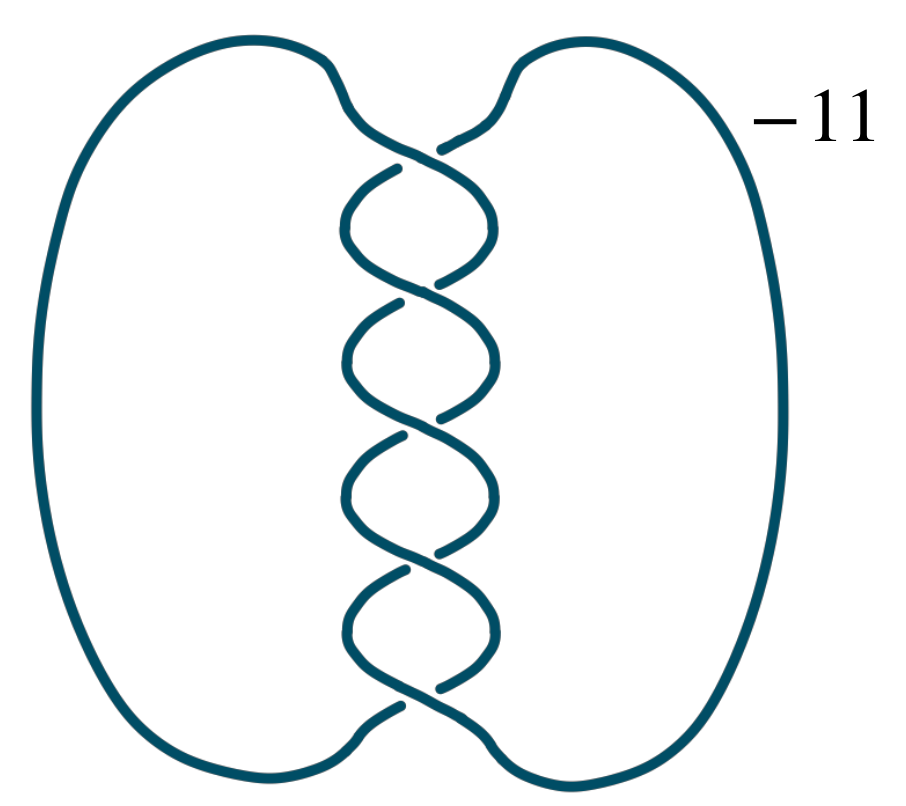}
\caption{Torus knot $T(-5,2)$.}
\label{t(5,2)}
\end{figure}

\begin{figure}[ht!]
\centering
\begin{subfigure}[t]{.45\textwidth}
  \centering
 \includegraphics[scale=0.45]{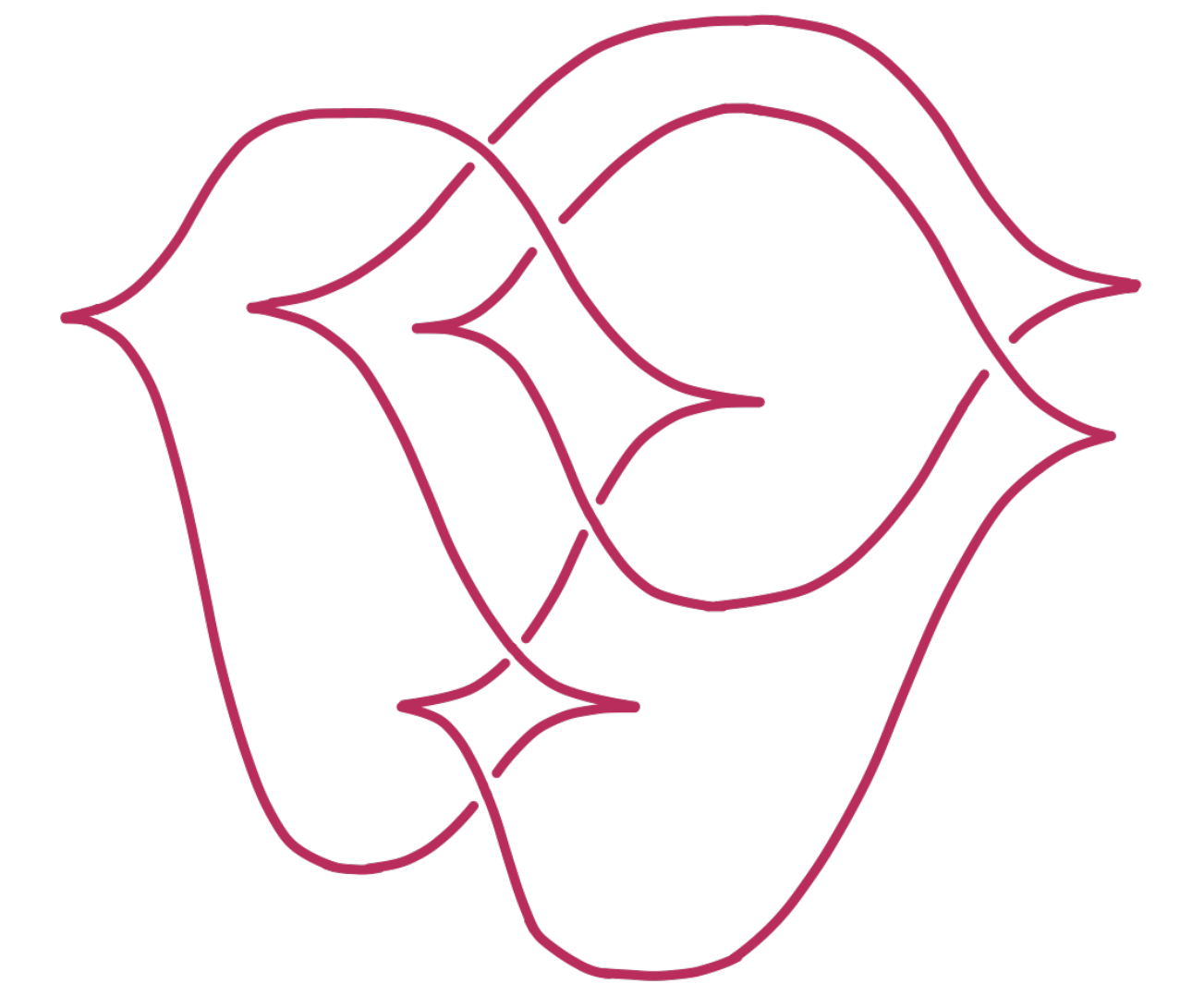}
\caption{Virtually overtwisted structure.}
  \label{t(5,2)1}
\end{subfigure}
\begin{subfigure}[t]{.45\textwidth}
  \centering
 \includegraphics[scale=0.5]{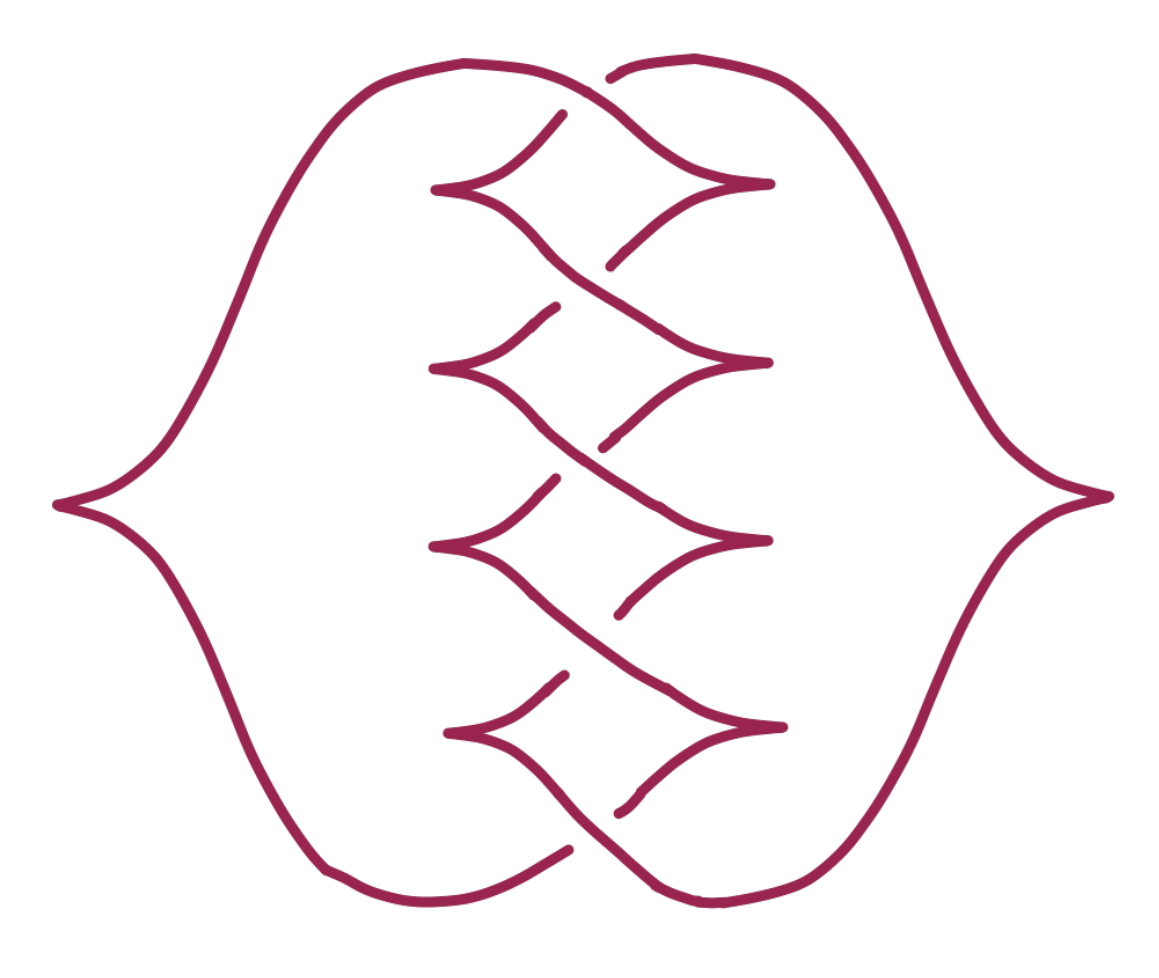}
  \caption{Universally tight structure.}
  \label{t(5,2)2}
\end{subfigure}
\caption{Different Legendrian representatives of $T(-5,2)$.}
\label{legrepr}
\end{figure}

In order to encode the Stein structure of $(X,J)$ in this handle decomposition we need a Legendrian representative of $T(-5,2)$ with Thurston-Bennequin number equal to $-10$: combining \cite[Theorem 4.3]{torusknotsimple} and \cite[Theorem 4.4]{torusknotsimple} we see that there are just two such Legendrian isotopy classes which maximize the Thurston-Bennequin number (equal to $-10$), distinguished by the rotation numbers, respectively $\pm 1$ and $\pm 3$ depending on the orientations (see Figure \ref{legrepr}). We want to understand which one suits to our case. Let $J_1$ and $J_2$ be the two Stein structures on $X$ described respectively by Figures \ref{t(5,2)1} and \ref{t(5,2)2}. 

\begin{prop} 
The Stein domain $(X,J)$ with a handle decomposition consisting of a single 2-handle is the one described by Figure \ref{t(5,2)1}, i.e. $J=J_1$. 
\end{prop}

\begin{prf} 

Call $\xi$ the contact structure described by Figure \ref{contactxi}. Remember that the two open book decompositions of Figures \ref{pagehopf(3,4)} and \ref{pagelantern} represent the same contact structure.
To prove the proposition, it is enough to check that the induced contact structure on $\partial (X,J_1)$ is isotopic to $\xi$. This is achieved by computing the 3-dimensional invariant:
\[d_3(\partial (X,J_1))=\frac{1}{4}(c_1(X,J_1)^2-3\sigma(X)-2\chi(X)).\]
If we call $K$ the Legendrian knot of Figure \ref{t(5,2)1}, then the first summand is given by 
\[\rot(K)\cdot [-11]^{-1}\cdot\rot(K),\]
while $\sigma(X)=-1$ and $\chi(X)=2$. By putting everything together we obtain
\[d_3(\partial (X,J_1))=-\frac{3}{11}.\]
On the other hand, the link of Figure \ref{hopf(3,4)} gives a Stein filling $(W,J_0)$ of $(L(11,4),\xi)$ with two 2-handles such that 
\[c_1(W,J_0)^2=[1,-2]\cdot Q_W^{-1}\cdot [1,-2]^T,  \qquad \sigma(W)=-2,	\qquad \chi(W)=3,\]
where $Q_W$ is the matrix of the intersection form, which is just the linking matrix
\[  \begin{bmatrix}
   -3 & 1 \\
   1 & -4 \\
   \end{bmatrix}.\] 
The computation shows again that
\[d_3(\partial (W,J_0))=-\frac{3}{11}.\]
Moreover, in the case when the rotation number of the second component of the Legendrian link is 0, the $d_3$ invariant of the resulting contact structure is $-1/11$. According to Honda's classification of tight contact structures on $L(11,4)$, this computation covers all the three (up to contactomorphism) possible cases, see next paragraph for computation in the universally tight case.

Therefore, we conclude that
\[(L(11,4),\xi)=\partial (W,J_0)=\partial(X,J_1).\]
Hence $J=J_1$, as wanted.
\end{prf}

\begin{figure}[ht!]
\centering
\includegraphics[scale=0.4]{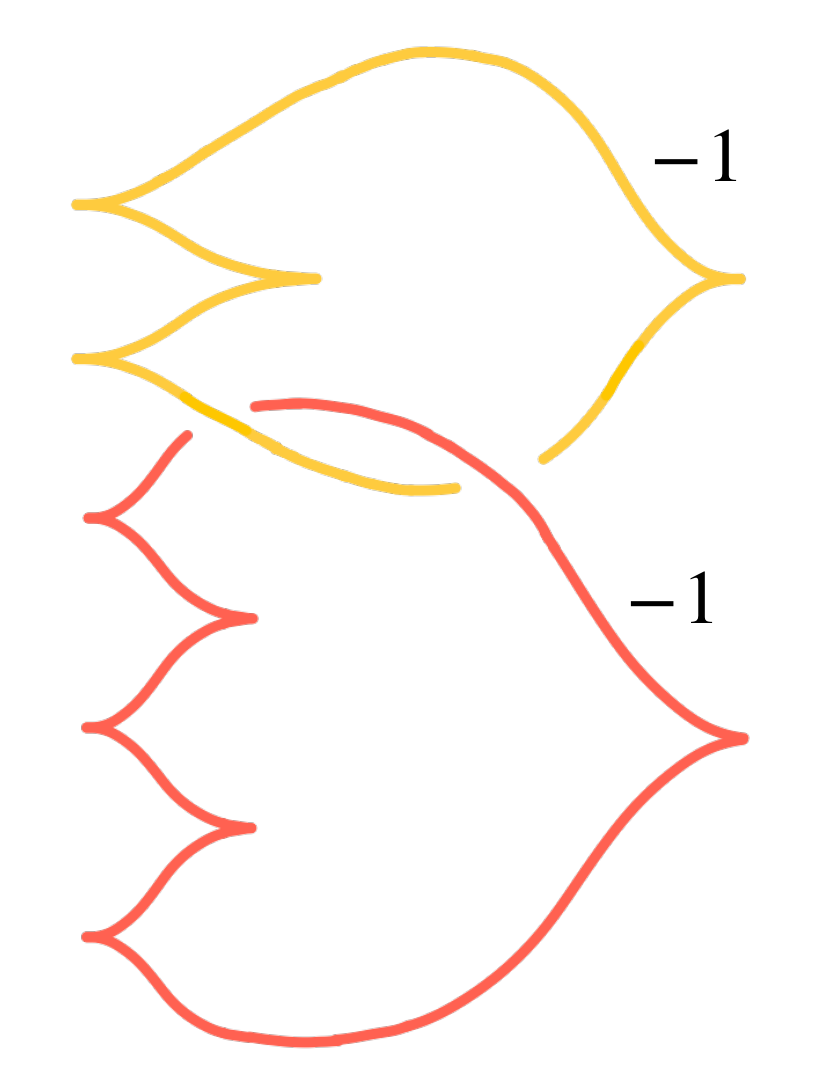}
\caption{Universally tight contact structure on $L(11,4)$.}
\label{hopf34ut}
\end{figure}

To conclude, we check that the other Legendrian representative of the torus knot $T(-5,2)$ with Thurston-Bennequin number $-10$ gives a different contact structure: by performing contact $(-1)$-surgery on the Legendrian knot of Figure \ref{t(5,2)2}, we get $L(11,4)$ with a universally tight structure. This is proved by comparing its $d_3$ invariant with the one computed from Figure \ref{hopf34ut}.

\noindent In both cases we get 
\[d_3=-\frac{5}{11}. \]
Therefore the two different Stein fillings of $(L(11,4),\xi_{ut})$ are described by the handle diagrams of Figures \ref{hopf34ut} and \ref{t(5,2)2}.

\bibliographystyle{alpha} %alpha, apalike

\bibliography{Fillings_hopf_link.bib}
\thispagestyle{plain}

\Addresses

\end{document}